%
\let\ReallyTheEnd=\end
\def\end{\vfill\eject}
%

\input fig.tex
\input psfig
\def\QP{\narrower\smallskip\noindent}
\def\ref{\hangindent=1pc \hangafter=1 \noindent}

\def\[{$\,}
\def\]{\,$}
\def\C{{\bf C}}
\def\R{{\bf R}}
\def\Z{{\bf Z}}

\def\={\;=\; }

\font\bit=cmssi12 at 12truept
\font\tenmsy=cmbsy10
\mathchardef\ssm="7872
\mathsurround = 1pt
\abovedisplayskip=6pt
\belowdisplayskip=6pt
\parskip=2pt

%

\input fig.tex
\input psfig
\def\QP{\narrower\smallskip\noindent}
\def\ref{\hangindent=1pc \hangafter=1 \noindent}

\def\[{$\,}
\def\]{\,$}
\def\C{{\bf C}}
\def\R{{\bf R}}
\def\Z{{\bf Z}}

\def\={\;=\; }

\font\bit=cmssi12 at 12truept
\font\tenmsy=cmbsy10
\mathchardef\ssm="7872
\mathsurround = 1pt
\abovedisplayskip=6pt
\belowdisplayskip=6pt
\parskip=2pt

\input psfig
\def\endproof{  \rlap{$\sqcup$}$\sqcap$ \smallskip}
\def\C{{\bf C}}
\def\QP{\medskip\leftskip=.4in\rightskip=.4in\noindent}
\def\[{$\,}
\def\]{\,$}
\def\Aut{ {\rm Aut}}
\def\crt{^{\cal C}}
\def\cl{{\cal C}}
\def\hl{{\cal H}}
\def\p{{\cal P}}
\font\bit=cmssi12 at 12truept
\font\figf=cmbx10 at 10truept


\centerline {\bf Hyperbolic components in spaces of polynomial maps.}
\smallskip

\centerline{J. Milnor, SUNY Stony Brook, February 1992}\smallskip

\centerline{With an Appendix by A. Poirier.}\smallskip

{\QP {\bf Abstract.} \bit We consider polynomial maps \[f:\C\to\C\] of
degree \[d\ge 2\],
or more generally polynomial maps from a finite union of copies of
\[\C\] to itself which have degree two or more on each copy.
In any space \[\p^{S}\] of suitably normalized maps of this type, the
post-critically bounded maps form a compact subset \[\cl^{S}\] called
the {~\bf connectedness locus}, and the hyperbolic maps in \[\cl^{S}\]
form an open set \[\hl^{S}\] called the {~\bf hyperbolic} connectedness
locus. The
various connected components \[H_\alpha\subset \hl^{S}\] are called
{~\bf hyperbolic components}. It is shown that
each hyperbolic component is a topological cell,
containing a unique post-critically finite map which is called its
{~\bf center point}. These hyperbolic components
can be separated into finitely many distinct ``types'',
each of which is characterized by a suitable
{~\bf reduced mapping schema}
\[\bar S(f)\]. This is a rather crude invariant, which depends only on
the topology of \[f\] restricted to the complement of the Julia set.
Any two components with the same reduced mapping schema are
canonically biholomorphic to each other. There are similar statements for
real polynomial maps, or for maps with marked critical points.
\medskip}

\centerline{\bf \S1. Introduction.}\smallskip

First some mild generalizations of standard concepts.
(See for example [D1], [DH1], [M2].)
Let \[M\] be the disjoint union of finitely many
copies of the complex numbers \[\C\],
and let \[f:M\to M\] be a proper holomorphic map, of degree at least two
on each component. The {\bit filled Julia set\/} \[K(f)\subset M\] is the
set of all points whose forward orbit is contained in some compact subset
of \[M\]. Such a map is {\bit post-critically bounded\/} if every critical
point is contained in \[K(f)\], and is {\bit
post-critically finite\/} if the orbit of every critical point is finite.
Define the {\bit fully invariant Julia}\footnote
{${}^1$}{This form of the
definition is actually due to Fatou rather than Julia. In
our mildly generalized context, note that
the closure of the repelling periodic set as considered by Julia
may well be strictly smaller than this fully invariant set \[J(f)\].}
{\bit set\/} \[J(f)\] to be
the complement of the region in \[M\]
where the collection of iterates of \[f\] forms a normal family.
Extending the classical theory of Fatou and Julia, it is easy
to check that \[J(f)\]
coincides with the topological boundary \[\partial K(f)\]. Furthermore,
the intersection of \[K(f)\] or \[J(f)\] with each component
of \[M\] is connected if and only if \[f\] is post-critically bounded.
Such a map is {\bit hyperbolic\/} (ie., hyperbolic on its
Julia set) if and only if every critical orbit converges to an attracting
cycle. Here a cycle (=periodic orbit) is called {\bit attracting\/}
if and only if its {\bit multiplier\/} \[\lambda\], the first
derivative around the orbit, satisfies \[|\lambda|<1\], so that it attracts
all orbits in a neighborhood.\smallskip

First consider the classical and most important case \[M=\C\].
A polynomial map \[f:\C\to\C\] of degree \[d\] is
{\bit monic\/} and {\bit centered\/} if it has the form
$$	f(z)~ =~ z^d + a_{d-2}z^{d-2} + \cdots + a_1z + a_0 \;.	$$
For \[d \ge 2\], let \[\p^{d-1}\] be the complex \[(d-1)$-dimensional affine
space consisting of\break all polynomial maps from \[\C\] to
itself which are monic and centered of degree \[d\].\break ({\bf Caution:} It might
seem more natural to index such a space by the degree \[d\]. However, for our
purposes it will be much more convenient
to index by the total number of critical points, which is \[d-1\].)
Note that every polynomial
map from \[\C\] to itself is conjugate under an affine change of variable
to a map \[f\] which is monic and centered. Furthermore, this \[f\] is uniquely
determined up to the action of the group \[G(d-1)\] of \[(d-1)$-st roots
of unity, where each \[\eta\in G(d-1)\]
acts on \[f\in \p^{d-1}\] by the transformation
$$	f(z) \mapsto f(\eta z)/\eta	\,.	$$
Compare 2.7 below.\smallskip

By definition, the {\bit connectedness locus\/} \[\cl^{d-1} \subset \p^{d-1}\] is the
compact set consisting of all polynomials \[f \in \p^{d-1}\] for which
\[K(f)\] is connected, or equivalently contains all critical points.
We define the
{\bit hyperbolic connectedness
locus\/} \[\hl^{d-1} \subset \cl^{d-1}\] to be the open set consisting
of all \[f\in \p^{d-1}\] for which the orbits of all critical points are
not only bounded, but actually converge to attracting periodic orbits.
\smallskip

{\bf Remarks.} The connectedness locus
\[\cl^1\] for quadratic maps, better known as the {\bit Mandelbrot set\/},
has been extensively studied by Douady and Hubbard.
The locus \[\hl^1\], consisting
of hyperbolic maps in \[\cl^1\], was considered somewhat earlier
by Brooks and Matelski. The {\bit cubic connectedness locus\/} \[\cl^2\]
has been studied by Branner and Hubbard.  In both the quadratic and the cubic
cases, an important result is that this connectedness
locus is {\bit cellular\/}, ie., equal to the intersection of a
strictly nested family of closed topological cells. The corresponding
statement for maps of higher degree has been proved by Lavaurs.
In all cases, it is conjectured that \[\hl^{d-1}\]
is everywhere dense in \[\cl^{d-1}\], and coincides with the
interior of \[\cl^{d-1}\]; however these statements have not been
proved even in the simplest case \[d=2\].\smallskip

Each connected component  \[H_\alpha\subset\hl^{d-1}\] is called a {\bit
hyperbolic component~} in \[\cl^{d-1}\]. For degree \[d=2\], each map
\[f\in\hl^1\] has a unique attracting periodic orbit. Let \[\lambda_f\] be its
multiplier. Douady and Hubbard [D1], [DH1]
have shown that every hyperbolic component \[H_\alpha\subset\hl^1\] is
canonically biholomorphic to the open unit disk \[D\] under the correspondence
$$	f\;\mapsto \;\lambda_f\in D\,.	$$
Rees [R1] has studied hyperbolic components for quadratic rational maps,
and McMullen [McM] has implicitly made a corresponding study for rational
maps of arbitrary degree with connected Julia set.
This paper is concerned with the special case of
polynomial maps.\smallskip

Section 2 will begin by defining the {\bit reduced mapping schema}
$$ \bar S(f)\;=\;(|\bar S|\,,\,\bar F\,,\,\bar w) $$
associated with any
hyperbolic polynomial map \[f\]. This consists of a finite set \[|\bar S|\]
having one point \[v\] for each component \[W_v\] of the interior of \[K(f)\]
which contains critical points, together with a ``first return map'' \[\bar F\]
from \[|\bar S|\] to itself, and an integer valued {\bit critical weight
function\/}, where \[\bar w(v)\ge 1\] is the number of
critical points in \[W_v\], counted with multiplicity.
%
To each reduced mapping schema \[ S\] there
is associated an affine space \[\p^{ S}\] of polynomial mappings whose
{\bit principal hyperbolic component} \[H_0^{ S}\] provides a standard model
for hyperbolic components with schema \[ S\]. This standard model has a finite
group \[\bar G(S)\cong\Aut(H_0^{ S})\] of automorphisms.\smallskip

Sections 3 and 4 show that each hyperbolic component \[H_\alpha\], in any space
\[\p^{S_1}\], is a topological cell, and that it
contains a unique post-critically
finite map, called its {\bit center\/} point
\[f_\alpha\]. (Compare [McM],
[R1].) The proof, using Douady-Hubbard surgery, shows
also that every hyperbolic component with reduced schema isomorphic to \[S\]
is canonically
diffeomorphic to a standard model \[B(S)\], which consists of certain
collections of Blaschke products. The diffeomorphism is unique up to
composition with an element of the group \[\bar G(S)\], which acts on
\[B(S)\]. In
particular, any two hyperbolic components with the same reduced schema
are canonically diffeomorphic to each other.  Section 5 sharpens this
statement by showing that they are canonically biholomorphic.
Section 6 discusses analogous results for polynomial mappings with
real coefficients, or more generally for ``real forms'' of complex polynomial
mappings. Section 7 studies polynomial mappings which have been
{\bit critically
marked\/} by specifying an ordered list of their critical points. It
is shown that all of the principal results carry over to the critically
marked case. The Appendix, by Alfredo
Poirier, shows that every possible reduced
mapping schema actually occurs as the schema associated with some critically
finite hyperbolic map \[f:\C\to\C\].
\smallskip

{\bf Remark.} Both
the statement that each hyperbolic component is a topological
cell, and the statement that it has a preferred center point,
are strongly dependent on the fact that we consider only maps with
connected Julia set. In the
case of quadratic rational maps, Rees shows that the unique
hyperbolic component consisting of maps with disconnected Julia set
has a more complicated topology.
(Compare [M4].) For polynomial
maps outside the connectedness locus, Blanchard, Devaney and Keen
describe a very rich topology within the open set consisting
of hyperbolic maps for which all critical orbits escape to infinity.\smallskip

The present work is a fairly straightforward extension of ideas originated
by Douady, Hubbard, McMullen, Rees and others. I am particularly
grateful to Branner and Douady for their
considerable help.\bigskip\vfil

\centerline{\bf \S2. The Mapping Schema of a Hyperbolic Component.}
\medskip

\noindent{\bf Definition 2.1.} By a {\bit mapping
schema} \[\;S\,=\,(|S|\,,\,F\,,\,w)\;\] we mean:

\ref (1) a finite set \[|S|\] of points, 
together with

\ref (2) a function \[F=F_S\] from \[|S|\] to itself, and also

\ref (3) a ``weight function'' \[w=w_S\] which assigns an
integer \[w(v) \ge 0\] called the {\bit critical weight\/}
to each \[v \in |S|\].

\noindent 
Equivalently, such
a mapping schema can be represented by a finite graph with one vertex for
each \[v \in |S|\], and with exactly one directed edge \[e_v\]
leading out from each
vertex \[v\] to a vertex \[F(v)\]. By definition, the {\bit degree\/}
associated with the edge \[e_v\] (or with the vertex
\[v\]) is the integer \[d(v)=w(v)+1\ge 1\].
The possibility that \[v=F(v)\], so that \[e_v\]
is a closed loop, is not excluded.
The vertex \[v\] is called
a ``critical point'' if \[w(v)\ge 1\], and a ``multiple critical point'' if
\[w(v)\ge 2\]. The sum
\[\;w(S)=\Sigma_{v\in |S|}\;w(v)\;\] is called the {\bit total critical
weight\/} of the schema \[S\].
\eject

Such a mapping schema is {\bit reduced\/} if every vertex is critical.
Suppose that we start with a
not necessarily reduced mapping schema \[ S\] which satisfies
the following very mild condition: {\it Every cycle in \[ S\]
must contain at least one critical point.} Then there is an {\bit
associated reduced
schema\/} \[\bar S\] which is obtained from \[ S\] simply by discarding
all vertices of weight zero and
shrinking every edge of degree one to a point. (Compare Figure 1.)
Note that \[ S\] and \[ \bar S\] have the
same critical weight. In the special case where \[ S\] is itself a
reduced mapping schema, evidently \[S= \bar S\].

\midinsert
\centerline{\psfig{figure=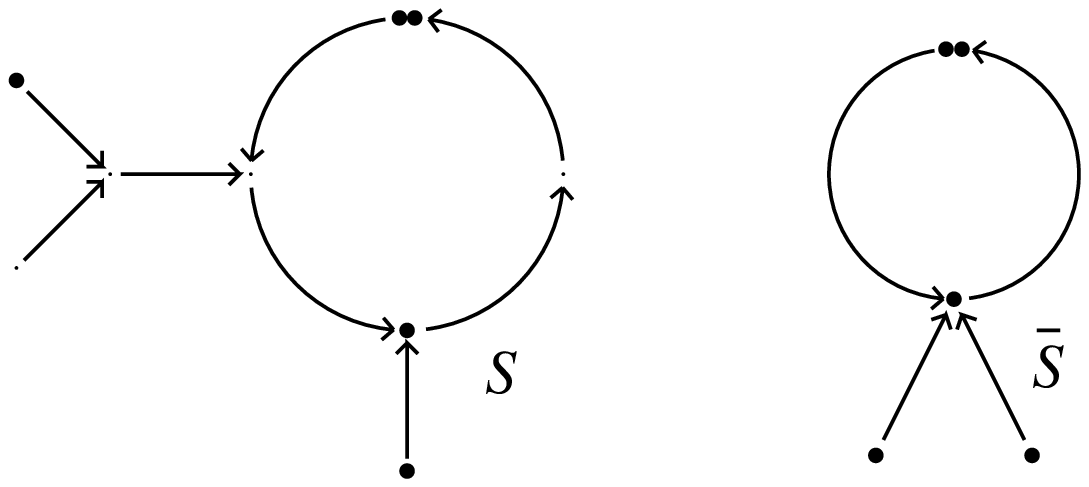,height=2in}}
{\QP\figf Figure 1. A mapping schema of total weight five, and the associated
reduced schema of weight five.
The heavy dots represent critical points, and the double
heavy dots represent critical points of weight two.\par}
\endinsert

Each mapping schema can be conveniently described by a symbol as follows.
The symbol \[(w)\]  stands for the graph with a single vertex of weight \[w\]
and a single edge which is a closed loop of degree \[w+1\],
while \[(w_1,\ldots,w_k)\], or any cyclic permutation of this symbol,
stands for a graph with \[k\]
vertices of weights \[w_1,\ldots,w_k\] arranged in a cycle. Note that every
connected mapping schema contains exactly one such cycle. To indicate
that additional vertices with weights \[w'_1,\ldots,w'_m\] map to the vertex
of weight \[w_i\], we insert the expression \[\{w'_1,\ldots,w'_m\}\]
immediately in front of the symbol \[w_i\]. Continuing this construction
inductively, we obtain an appropriate symbol for any connected mapping schema. 
As an example, the
schema \[S\] of Figure 1 has symbol \[\;\big(\{\{1,0\}0\}0\,,\,\{1\}1\,,\,
0\,,\,2\big)\;\] up to cyclic permutation,
and the associated reduced schema \[\bar S\] has symbol \[(\{1,1\}1\,,\,2)\].
Finally,
using the notation \[S+S'\] for a disjoint union, we obtain symbols for
disconnected schemata also.\smallskip

\midinsert
\centerline{\psfig{figure=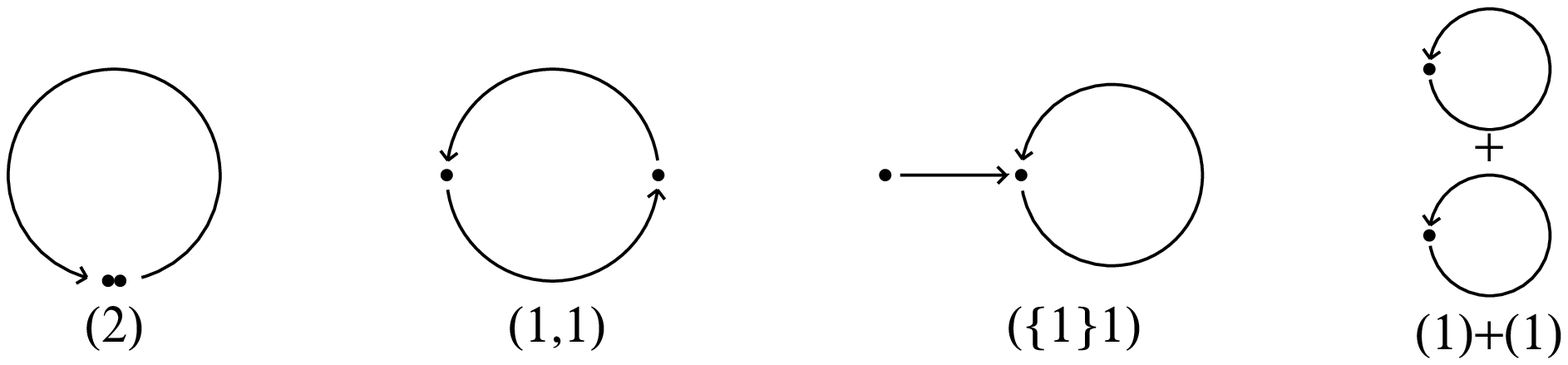,height=1.1in}}
\medskip
\centerline{\figf Figure 2. The reduced schemata of weight two.}
\endinsert

Up to isomorphism, there is just one reduced mapping schema of total weight
one, with symbol \[(1)\]; while there are four reduced schemata
\[(2)\,,\;(1,1)\,,\;(\{1\}1)\], and \[(1)+(1)\]
of weight two. (Compare [R1]. In [M3] the letters
\[A,\,B,\,C,\,D\] were used for these four types.) I don't know any
general formula for the number \[N(w)\] of distinct reduced schemata of
weight \[w\]. However, for \[n\le 5\] this number can be computed as follows:
\smallskip

\hskip80pt $N(1)\;=\;\;1$

\hskip80pt $N(2)\;=\;\;4\;\;=\;\;3\;+\;1$

\hskip80pt $N(3)\;=\;12\;=\;\;8\;+\;3\;+\;1$

\hskip80pt $N(4)\;=\;42\;=\;24+14\,+\;3+1$

\hskip80pt $N(5)\;=138=\;72+48+14+3+1\;\;.$
\smallskip
\noindent Here the first summand gives the number of connected schemata, the
next gives the number with two components, and so on. This can be proved simply by
constructing an exhaustive list.\smallskip

Any hyperbolic polynomial map  \[f\] from \[\C\] to \[\C\], or from a finite
union of copies of \[\C\] to itself of degree two or more on each copy,
gives rise to an {\bit associated full mapping schema\/}
\[ S=( |S|\,,\, F\,,\, d)\] by the following construction.
Let \[W(f)\] be the union of the basins
of attraction for all attracting periodic orbits of \[f\] in \[\C\].
(Equivalently, since \[f\] is hyperbolic,
\[W(f)\] is the interior of the filled Julia set.) Thus \[W(f)\] is the
bounded open set consisting of all points \[z\] whose forward orbit
converges to a attracting periodic orbit. It is not difficult to check that:
\smallskip

(a) {\it each connected component of \[W(f)\] is simply connected, and hence is
conformally isomorphic to the open unit disk; and}\smallskip

(b) {\it the map \[f\] carries each component \[W_\alpha\] of \[W(f)\] properly
onto some component \[f(W_\alpha)=W_\beta\] of \[W(f)\] by a map of degree \[d_\alpha
\ge 1\].}\smallskip

Furthermore, \[f\] is bijective on a given component \[W_\alpha\] if and
only if \[d_\alpha=1\], or in other words if and
only if there are no critical points in \[W_\alpha\]. Let \[W^{\cal
PC}\] be the union of those components \[W_v\] of \[W(f)\] which
contain {\bit post-critical\/} points of \[f\], that is points which are either
critical or belong to the forward orbit of some critical point.
 \smallskip

{\bf Definition 2.2: The full mapping schema of \[f\].} This schema
\[ S(f)\] has one
vertex \[v\] corresponding to each post-critical component \[W_v
\subset W^{\cal
PC}$. The {\bit weight} \[ w(v)\] is defined to be the number
of critical points in \[W_v\], counted with multiplicity. Each vertex \[v\]
is joined to a vertex \[F(v)\] by an edge \[e_v\] of degree \[w(v)+1\],
where \[F(v)\] is the vertex associated with the component \[f(W_v)=W_{F(v)}
\subset W^{\cal PC}\].
\smallskip

{\bf Definition 2.3: The reduced schema of \[f\].} In practice, we will
always work with the
associated reduced mapping schema \[\bar S(f)\], which can be
constructed from \[ S(f)\] as above, or can be
constructed directly as follows. {\it Let \[ W\crt\subset W^{\cal PC}\] be the
union of the {\bit critical components\/} of \[W(f)\], that is the union of
those components \[W_v\] which
contain critical points of \[f\], and let \[f\crt : W\crt  \to W\crt \] be
the first return map.} In other words, \[f\crt (z)\] is equal to the first
of the points
$$	 f(z)\,,\;\; f^{\circ 2}(z)=f(f(z))\,,\;\;f^{\circ 3}(z)= 
	f(f(f(z)))\,,\;\; \cdots	$$
which belongs to a critical component of \[W(f)\]. It follows from
the classical theory of Fatou and Julia that every attracting periodic
orbit contains at least
one point of \[W\crt \], so that \[f\crt\] is defined. By definition,
the schema \[\bar S(f)\] has one vertex
\[v\] for each critical component \[W_v\subset W\crt \], with a
directed edge leading from each vertex \[v\] to the vertex \[\bar F(v)\] where
the first return map \[f\crt\] on \[W\crt\] carries \[W_v\] onto \[W_{\bar
F(v)}\]. The weight \[w(v))\ge 1\] is equal to the number of critical
points of either
\[f\crt\] or \[f\] in \[W_v\], counted with multiplicity.\smallskip

Thus the total
critical weight \[w\] of \[ S(f)\] or of \[\bar S(f)\] is equal
to the total
number of critical points of \[f\], counted with multiplicity.
We will see in the Appendix that all possible reduced schemata
arise in this way: 
{\it Given any reduced mapping schema \[S_0\] with total
critical weight \[w\] there exists a polynomial map \[f:\C\to\C\]
of degree \[d=w+1\] which is hyperbolic, with reduced
schema \[\bar S(f)\] isomorphic to \[S_0\].} 
We will say briefly that such an \[f\] has {\bit type\/} \[S_0\].

If we allow
mappings from a finite union of copies of \[\C\] to itself, then such
a map \[f\] can be constructed trivially as follows.
\smallskip

{\bf Definition 2.4. The universal polynomial model space.}
To each reduced mapping schema \[S=(|S|\,,\,F\,,\,w)\] there
is associated a complex affine space \[\p^S\]
of polynomial maps as follows. Form the disjoint union \[|S|\times\C\]
of \[n\] copies
of the complex numbers \[\C\], where \[n\] is the number of vertices of
\[S\]. In other words, replace each vertex of \[S\] by a copy of \[\C\].
Let \[\p^S\] be the space consisting of all maps \[f\] from \[|S|\times\C\]
to itself such that the restriction of \[f\] to each component \[v \times
 \C\] is a monic centered polynomial map of degree \[\;d(v)=w(v)+1\,\],
taking values in \[F(v)\times\C\].
Note that the complex dimension of this affine space
\[\p^S\] is equal to the total critical weight \[w(S)\].

{\bf Remark.} The principle that the dimension of a suitably normalized
space of holomorphic maps is equal to the (generic) number of critical points
seems to be true in a number of interesting
cases. Here are two further examples. For polynomials
of degree \[d\ge 2\] with just one multiple critical point, there is a
normal form \[z\mapsto z^d+c\] depending on one parameter. On the
other hand, the space of all holomorphic conjugacy classes of
rational maps \[f:\C\cup\infty\to\C\cup\infty\]
of degree \[d\ge 2\] has complex dimension \[2d-2\], and each such map
has exactly \[2d-2\] critical points counted with multiplicity.

In the special case where \[S\] is the schema \[(w)=(d-1)\] with just one vertex,
clearly the space \[\p^{(w)}\] of 2.4
coincides with the space \[\p^{d-1}\] described in \S1. For
any \[S\], the connectedness locus \[\cl^S \subset \p^S\] and the hyperbolic
set \[\hl^S \subset \cl^S\] can be defined just as in the special case
\[S=(w)\]. By definition, a map \[f\in \p^S\] belongs to the connectedness
locus \[\cl^S\] if and only if its filled Julia set
\[K(f)\subset |S|\times\C\] intersects each component
\[v\times\C\] in a connected set, or if and only if every critical point
has bounded orbit. Each hyperbolic map \[f\in \hl^S\] will
have its own reduced schema \[\bar S(f)\]. In general, \[\bar
S(f)\] will not be
isomorphic to the given schema \[S\], although it will have the same
total critical weight. (Compare 2.11.)\smallskip

There is a preferred base point \[f_0\] in \[\hl^S\], namely
the unique map
$$	f_0(v\,,\,z)~ =~ (F(v)\,,\,z^{d(v)})	$$
whose constituent polynomials are all monic monomials. Clearly the
mapping schema \[S(f_0)=\bar S(f_0)\] associated with this central point
can be identified with
the given reduced schema \[S\]. The hyperbolic component \[H_0^S\]
containing \[f_0\] will be called the {\bit principal\break
hyperbolic component}
of \[\cl^S\]. We will use \[H_0^S\] as a standard model for
hyperbolic components with reduced schema isomorphic to \[S\].
\smallskip

{\bf Definition 2.5. The symmetry groups \[G(S)\to \bar G(S)\].} By the
group \[\Aut(S)\] of {\bit automorphisms\/} of the schema \[S=(|S|\,,\,F\,,\,w)\],
we mean the finite group consisting
of all permutations of the set \[|S|\] of vertices which commute
with the map \[F\] and preserve the weight \[w(v)\]. As examples, for the
graphs shown in Figures 1b, 2b, and 2d, this automorphism group
is cyclic of order two, but for the remaining graphs in Figures 1 and 2 it
is trivial. (Compare 2.10 below.) Each such
automorphism gives rise to an automorphism of the spaces
$$	H_0^S \;\subset\; \hl^S\;\subset\;\cl^S \;\subset\; \p^S	$$
which is {\bit linear\/} in the sense that it preserves the affine structure
and fixes the base point \[f_0\]. However, these spaces
have a possibly larger group of linear automorphisms, as follows.
We first introduce the {\bit full symmetry group\/} \[G(S)\]
consisting of all maps \[g\] from \[|S|\times\C\] to itself which are
bijective, commute with \[f_0\], and carry each component
\[v\times\C\] linearly onto some \[v'\times\C\].

{\QP {\bf Lemma 2.6.} \it Conjugation by an element \[g\in G(S)\] carries
each polynomial map \[f\in \p^S\] to a map \[g\circ f\circ g^{-1}\in \p^S\].
The resulting action of the group \[G(S)\] on the space \[\p^S\] is linear, and
preserves the subsets \[H_0^S\subset \hl^S\subset \cl^S\].\medskip}

The proof will be given below.

{\bf Caution.} This action need not be effective:
Some elements of \[G(S)\] may commute with all
maps in \[\p^S\], and hence act trivially on the space \[\p^S\].
We will use the notation \[\bar G(S)\] for the {\bit effective symmetry
group\/}, that is the quotient group of \[G(S)\] which operates
effectively on the spaces \[\hl_0^S\subset\hl^S\subset\cl^S\subset\p^S\].
For details, see 2.9 below.
\smallskip

In order to justify the definition of \[\p^S\] and \[G(S)\], let us
note that any possible dynamical behavior which can
be obtained by nowhere linear polynomial maps from a finite union
of copies of \[\C\] to itself can already be realized by some mapping \[f\]
which belongs to an appropriate affine space \[\p^S\] of monic centered
polynomials. (This is an immediate generalization of the statement that
every polynomial map \[\C\to \C\] of degree \[d\ge 2\] is affinely conjugate
to a monic centered map.)

{\QP {\bf Lemma 2.7.} {\it Let \[M\] be a finite union of disjoint
copies of \[\C\], and let
$$	\phi : M\; \to\; M	$$
be a map
whose restriction to each copy of \[\C\] is given by a polynomial of\break
degree \[\ge 2\]. Then there is a mapping schema \[S\], unique up to
isomorphism, and a bijection
\[h : M \to |S|\times\C\] which is
affine on each copy of \[\C\], so that the conjugate map
\[f=h \circ \phi \circ h^{-1}\]
belongs to the space
\[\p^S\]. Furthermore, \[h\] is unique up to composition with
elements of the finite group \[G(S)\] which acts on \[|S|\times\C\],
and \[f\] is unique up to the action of \[\bar G(S)\] on \[\p^S\].} \smallskip}

\noindent The proof is easily supplied. \endproof\smallskip

To describe the group \[G(S)\] more explicitly, let \[m\] be
the number of connected components of \[S\], and let \[\delta_k\] be the
product of the degrees around the unique cycle which is contained in the
\[k$-th component of \[S\]. Let \[\delta'\] be the product of all of the
remaining degrees, which do not lie in cycles.

{\QP{\bf Lemma 2.8.} \it
The group \[G(S)\] is finite, and splits canonically as a
semi-direct product
$$	1\to N(S)\to G(S)
	\buildrel{\displaystyle\leftarrow}\over\to \Aut(S)\to 1\,.	$$
Here \[\Aut(S)\] is defined to be the group of automorphisms of
\[S\] $($see $2.5)$, while \[N(S)\subset G(S)\] is
a normal abelian subgroup of order
\[	(\delta_1-1)\cdots(\delta_m-1)\,\delta' \,.	\]
This subgroup can be identified with the group of all functions \[\eta\]
from \[|S|\] to the unit circle which satisfy the identity \[\;\eta(F(v)) =
 \eta(v)^{d(v)}\].\smallskip}\smallskip

{\bf Proof of 2.8 and 2.6.} Let \[N(S)\] be the subgroup consisting of those elements
of \[G(S)\] which carry each copy \[v\times\C\] to itself. Then any element
\[g \in N(S)\] must have the form \[g(v,z)=(v\,,\,\eta(v)z)\], where the \[\eta(v)\]
are non-zero complex numbers. In fact, a brief computation shows that such a
map \[g\] commutes with \[f_0\]
if and only if the numbers \[\eta(v)\] satisfy the
identities \[\eta(F(v))= \eta(v)^{d(v)}\]. If \[v\] lies in the \[k$-th cycle of
\[S\], then it follows easily that \[\eta(v)\] must be a \[(\delta_k-1)$-st
root of unity.
Evidently, the correspondence \[f \mapsto g\circ f \circ g^{-1}\]
defines a linear action of the finite group \[N(S)\] on the affine space
\[\p^S\], preserving the subsets \[H_0^S\], \[\hl^S\], and \[\cl^S\].
Further details will be left to the reader. \endproof

{\bf Remark 2.9.} A similar argument shows that the {\bit effective symmetry
group\/} \[\bar G(S)\] can be identified with the quotient
\[G(S)/N_0(S)\], where \[N_0\subset N\]
is a direct sum of cyclic groups of order two,
consisting of all functions from \[|S|\] to
\[\{\pm 1\}\] which are trivial on the image \[F(|S|)\], and also trivial on
all vertices which have weight \[w(v)>1\]. In other words, an automorphism
\[g\in G(S)\] of the space \[|S|\times\C\] commutes with every map
\[f\in\p^S\] if and only if it belongs to this subgroup \[N_0(S)\]. {\it
It follows that \[G(S)\ne \bar G(S)\] if and
only if there exists a vertex of weight \[w=1\] which does not belong to
the image \[F(|S|)\subset|S|\].}\smallskip

{\bf 2.10. Examples.}
In the quadratic case \[S=(1)\], the symmetry group \[G(1)\] is
trivial. For the four schemata of weight 2, these groups can be tabulated as
follows,\medskip

\centerline{\vbox{\hrule\hbox{\vrule\vbox{
\halign{~\hfil #\hfil &\quad\hfil #\hfil &\hfil #\hfil &\hfil #\hfil &
\quad\hfil #\hfil~\cr
\noalign{\smallskip}
\[S\] & \[N(S)\] & G(S) & \[\Aut(S)\] & \[\bar G(S)\]\cr
\noalign{\smallskip\hrule\smallskip}
\[(2)\] & \[\Z/2\] & \[\Z/2\] & 0 & \[\Z/2\]\cr
\[(1,1)\] & \[\Z/3\] & order 6 & \[\Z/2\] & order 6\cr
\[(\{1\}1)\] & \[\Z/2\] & \[\Z/2\] & 0 & 0\cr
\[(1)+(1)\] & 0 & \[\Z/2\] & \[\Z/2\] & \[\Z/2\]\cr
\noalign{\smallskip}
}
}\vrule}\hrule}} \medskip

\noindent where \[G(1,1)=\bar G(1,1)\] is the dihedral group of order \[6\].
Here are four analogous examples, where the total weight \[w\] can be
arbitrary: For the schema
\[S=(w)\] with a single vertex, the group \[N=G=\bar G\] is cyclic of
order \[w\] and \[\Aut(S)\] is trivial. For the cycle \[S=(1,1,\ldots,1)\]
of total weight \[w\], the group \[\Aut(S)\] is cyclic of order \[w\]
and \[N(S)\] is cyclic of order \[2^w-1\], hence the group
\[G(S)=\bar G(S)\] is non-abelian of order \[w\,(2^w-1)\].
For the schema \[\;S=(\{\{\cdots\{\{\,1\}\,1\}\cdots \,1\}\,1\}\,1)\;\] with
vertices \[\;v_1\to v_2\to
\cdots\to v_{w-1}\to v_w$\raise 2.9pt\hbox{$\leftarrow$}$\!\!\!\!\supset\;,\;\]
the group \[\Aut(S)\] is trivial and \[N=G\] is
cyclic of order \[2^{w-1}\], while \[\bar G\] is cyclic of order \[2^{w-2}\].
Finally, for the sum \[(1)+\cdots+(1)\] of total
weight \[w\], the group \[G=\bar G=\Aut(S)\] is the
symmetric group of order \[w!\] and \[N\] is trivial. The proofs are not
difficult.\medskip

{\bf Remark 2.11.} An interesting relation between reduced mapping
schemata of the
same weight can be defined as follows. {\it Let us say that \[S \succ S'\] if and
only if the connectedness locus \[\cl^S\] contains a hyperbolic component
with reduced schema isomorphic to
\[S'\].} This is clearly a reflexive relation; that is \[S \succ S\] in
all cases.
\smallskip

As an example, for the four reduced schemata of critical weight \[w=2\],
it is not difficult to check that this partial ordering is transitive,
and can be described as follows.\break
We have
$$	(2)\; \succ\; (1,1)\; \succ\;\left\{\eqalign{&\;\; (\{1\}1)	\cr
	& (1)+(1)\;\;, \cr}\right.	$$
but no other relations, except as implied by reflexivity and
transitivity. (There is an analogous partial ordering
for the various real forms of these hyperbolic components. Compare 6.2, and
Figures 3--5.)
\smallskip

If \[S\succ S'\], where \[S'=(|S'|\,,\,F'\,,\,w')\] and
\[S=(|S|\,,\,F\,,\,w)\], then evidently there is an associated map \[\psi:
|S'|\to |S|\] which:

(1) preserves order, in the sense that if \[F'(v_1')=v_2'\] then some iterate
of \[F\] maps \[\psi(v_1')\] to \[\psi(v_2')\], and

(2) preserves weight, in the sense that
\[	\sum_{\psi(v')=v}\,w'(v')\;=\;w(v)	\]
for each \[v\in|S|\].

\noindent It is conjectured that \[S\succ S'\] if and only if there exists
such a map \[\psi\].
In particular, it is conjectured that this relation
is transitive. In many cases, a hyperbolic
component \[H_\alpha\subset \hl^S\] of type \[S'\] seems to be contained in
a complete copy of the connectedness
locus \[\cl^{S'}\], which is homeomorphically
embedded into \[\cl^S\] in such a way that the principal component
\[H_0^{S'}\] corresponds to \[H_\alpha\]. (This is a generalization of
the Douady-Hubbard concept of ``modulation'' or
``tuning'' for quadratic maps.) Whenever this happens, it follows 
that any type of hyperbolic component which occurs in \[\cl^{S'}\] must
also occur in \[\cl^S\]. Thus, in this special case,
the relations \[S\succ S'\] and
\[S'\succ S''\] certainly imply that \[S\succ S''\]. 
\smallskip

\smallskip

{\bf Caution.} It may happen that \[S \succ S'\] and \[S' \succ S\], even
though \[S\] is not isomorphic to \[S'\]. For example,
this is true for \[S=(\{2\}1,\{1\}1,1)\] and \[S'=(\{1\}1,
 \{2\}1,1)\]. The proof of this statement is an easy application of the
Realization Theorem for ``Hubbard forests'', as proved by Poirier [P3].

\end

\def\QP{\medskip\leftskip=.4in\rightskip=.4in\noindent}
\def\ref{\hangindent=1pc \hangafter=1 \noindent}
\def\cl{{\cal C}}
\def\hl{{\cal H}}
\def\p{{\cal P}}
\def\s{{\cal S}}
\def\endproof{  \rlap{$\sqcup$}$\sqcap$ \smallskip}
\def\C{{\bf C}}
\def\R{{\bf R}}
\def\Z{{\bf Z}}
\def\[{$\,}
\def\]{\,$}
\def\Aut{{\rm Aut}}
\def\crit{^{ C}}
\font\bit=cmssi12 at 12truept
\font\tenmsy=cmbsy10
\textfont8=\tenmsy
\mathchardef\ssm="7872
\pageno=10

\centerline{\bf \S3. Blaschke Products and the Model Space \[B(S)\].}
\bigskip

Henceforth, all mapping schemata will be reduced.
This section will describe a topological model, based on
Blaschke products, for hyperbolic components with mapping schema \[S\].
First recall some standard facts.\smallskip

Let \[D\] be the open unit disk in \[\C\]. For any \[a \in D\],
there is one and only one conformal automorphism \[\mu_a:\bar
D\to \bar D\] which
maps \[a\] to zero and fixes the boundary point 1, given by
$$      \mu_a(z) \;=\; k\,{z-a\over 1-\bar a z}
	\qquad{\rm with}\qquad k={1-\bar a\over 1-a}\,.	$$
It is easy to check that the proper holomorphic
maps from \[D\] onto itself are precisely the finite products of the form
$$	\beta(z)\; =\; c\,\mu_{a_1}(z)\,\cdots\,\mu_{a_d}(z)\,,	 $$
with \[|c|=1\].
Here \[d\] is the degree, \[\{a_1\,\ldots,\,a_d\}\] are the pre-images
of zero, and \[c=\beta(1)\]. Evidently every such map \[\beta\] extends
uniquely as a rational map from
the Riemann sphere \[\hat\C=\C \cup \infty\] onto itself. In particular,
\[\beta\] extends continuously over the
boundary circle \[\partial D\]. Note that
the extension of \[\beta\] over \[ \hat\C\] commutes with the inversion
\[z \mapsto 1/\bar z\]. In particular, if \[z_0\] is a fixed point or a
critical point of \[\beta\], then \[1/\bar z_0\] is also. 

{\QP{\bf Lemma 3.1.} \it
A proper map \[\beta\] of degree \[d \ge 2\] from the unit disk onto
itself induces a \[d$-to-one covering map from the circle \[\partial
D\] onto itself. Such a map \[\beta\] has at most one
fixed point in the open disk \[D\]. If there is an
interior fixed point, then the induced map on
\[\partial D\] is topologically conjugate to the linear map \[t\mapsto
td\] of the circle \[\R/\Z\]. In particular,
in this case there are exactly \[d-1\] boundary fixed points.\medskip}

(On the other hand, if there is no interior fixed point, then there must
be \[d+1\] boundary fixed points, counted with multiplicity.)\smallskip

{\bf Proof Outline.} For \[|z|=1\], a brief
computation shows that the logarithmic\break derivative
$$	d\log\beta(z)/d\log(z)\; =\; z\,\beta'/\beta	$$
is a sum of \[d\] terms, each of which
is real and strictly positive. It follows easily that \[\beta\]\break induces a
covering map on the boundary, and that any
equation of the form\break \[\beta(z)= {\rm constant}\in \partial D\;\;\] has
exactly \[d\] distinct solutions. If there is an interior fixed point,
then after conjugating by a conformal automorphism, we may assume that
the fixed point is \[z=0\], so that
$$      \beta(z)\; =\; c\,z\,\mu_{a_2}(z)\cdots\mu_{a_{d}}(z) \,.	$$
It follows that \[|\beta'(0)|=|a_2 \cdots a_{d}|<1\] and that
 \[|\beta(z)|<|z|\] for all \[z \ne 0\] in the open disk. Thus this fixed
point is attracting, and is unique within \[D\].
When there is such an interior fixed point,
it is not difficult to check that the logarithmic derivative
satisfies \[z\beta'(z)/\beta(z) > 1\] at all points of
\[\partial D\]. In other words, the map of \[\partial D\] onto itself
is strictly expanding.
In particular, every boundary fixed point \[z_0\] is
repelling, with \[\beta'(z_0)>1\]. Since the algebraic number of
fixed points in \[\hat\C\] is equal to \[d+1\], it then follows that there are
exactly \[d-1\] distinct fixed points on \[\partial D\].

More explicitly, we can set up a coordinate system on \[\partial D\] as
follows. Choose one of these boundary fixed points \[z_0\] and assign it the
coordinate \[t(z_0)=0\]. Now assign the \[d\] preimages of \[z_0\] under
\[\beta\] the coordinates \[t(z)=j/d\], numbering in counterclockwise order
around the circle from \[z_0\]. Similarly, assign the \[d^2\] preimages under
\[\beta^{\circ 2}\] the coordinates \[j/d^2\], and so on. Since \[\beta\]
is expanding on \[\partial D\], the iterated preimages of \[z_0\] are
everywhere dense on \[\partial D\], so this construction converges to
a well defined homeomorphism\break
\[t:\partial D\buildrel\approx\over\to\R/\Z\].
\endproof

Closely related is the statement that \[\beta\]
gives rise to a preferred \[\beta$-invariant
measure \[\ell=\ell_\beta\]
on \[\partial D\], with the properties that \[\ell(\partial D)=1\] and that
\[\beta\] maps any small interval of measure \[\ell(I)=\epsilon\]
to an interval of measure \[\ell(\beta(I))=\epsilon\, d\]. For example
\[\ell\] can be defined by the
formula \[\ell(I)=\lim_{k\to\infty}\;N(k)/d^k\], where \[N(k)\] is the
number of solutions to the equation \[\beta^{\circ k}(u)=1\] in the interval
\[I\subset \partial D\]. Note that this measure is ``balanced'' in the
sense that each of the \[d\] components of \[\beta^{-1}(I)\] has
measure \[\ell(I)/d\]. Evidently \[\ell\] corresponds to the
Lebesgue measure on \[\R/\Z\] under the conjugating homeomorphism.\break
(This \[\ell\] can also be characterized
as the unique invariant measure of maximal entropy.)\medskip

{\bf Definition.} Let \[\beta:D\to D\] be a proper holomorphic map of
degree \[d\ge 2\]. We will say that \[\beta\] is {\bit boundary-rooted~}
if its extension over \[\partial D\] fixes the point \[+1\], and
{\bit fixed point centered~} if \[\beta(0)=0\].

{\QP{\bf Lemma 3.2.} \it If \[\phi\] is a proper holomorphic map of degree
\[d\ge 2\] from \[D\] to itself, with a fixed point in the open disk \[D\],
then \[\phi\] is holomorphically conjugate to a map
\[h^{-1}\circ\phi\circ h\] which is boundary-rooted and fixed point centered.
In fact there are exactly \[d-1\] distinct choices for the conjugating
M\"obius automorphism \[h\]. The space \[B(d-1)\] consisting of all proper
holomorphic \[\beta:D\to D\]
which are boundary-rooted and fixed point centered of degree \[d\]
is a topological cell of dimension \[2d-2\].\medskip}

{\bf Proof.} Evidently \[h\] can be any M\"obius automorphism of \[\bar D\] which
maps \[0\] to the unique interior fixed point of \[\phi\], and maps \[+1\]
to one of the \[d-1\] boundary fixed points. 
Let \[\s_n(\C)\] be the \[n$-{\bit fold symmetric product\/},
consisting of unordered \[n$-tuples\break
\[\{a_1\,,\,\ldots\,,\,a_n\}\] of complex
numbers. This can be identified with the complex affine space consisting of
all monic polynomials of degree \[n\], under the correspondence
$$      \{a_1\, ,\,\ldots\,,\,a_n\} \mapsto (z-a_1)\cdots (z-a_n)= z^n
        -\sigma_1 z^{n-1} +-\cdots \pm\sigma_n \,, $$
where the \[\sigma_j\] are the elementary symmetric functions of \[\{a_1\,,
 \,\ldots\,,\,a_n\}\]. Thus \[\s_n(\C)\] is homeomorphic to \[\C^n\cong
\R^{2n}\]. Since \[\C\] is homeomorphic to the 2-cell \[D\],
it follows that \[\s_n(D)\] is also homeomorphic to \[\R^{2n}\].
Now consider the space \[B(w)\] consisting of all boundary-rooted
Blaschke products \[\beta\] of degree \[d=w+1\] which fix the origin. We
can write
$$     \beta(z) \;=\; z\, \mu_{a_1}(z)\cdots \mu_{a_w}(z) \,.        $$
Evidently this space is homeomorphic to the symmetric product \[\s_{w}(D)\],
and hence is a topological cell, homeomorphic to \[\R^{2w}\].\endproof

We will also need to find a normal form for Blaschke products under
one-sided composition with a M\"obius transformation. {\bf Definition:}
A proper holomorphic map \[D\to D\] of degree \[d\ge 2\] is
{\bit critically centered~} if the sum of its \[\;d-1\;\]
(not necessarily distinct) critical points is equal to zero.

{\QP{\bf Lemma 3.3.} \it If \[\phi:D\to D\] is proper and holomorphic
of degree \[d\ge 2\], then there exist exactly \[d\] distinct M\"obius
automorphisms \[h:D\to D\] for which the composition \[\beta'=\phi\circ h\] is
boundary-rooted and critically centered. The space \[B'\] consisting of all
boundary-rooted, critically centered proper holomorphic maps of degree
\[d\] is a topological cell of dimension \[2d-2\].\medskip}

\noindent
The proof will be based on the concept of ``{\bit conformal barycenter\/}''
for a collection of points in the disk:
 
{\QP {\bf Lemma 3.4.} {\it Given points \[c_1\,,\,\ldots\,,\,c_n\] in
a Riemann surface \[W\] isomorphic to \[D\],
there exists a conformal isomorphism \[\eta: W \to D\], unique up to a
rotation of \[D\], which takes the 
\[c_j\] to points \[\eta(c_j)\] with sum \[\eta(c_1)+\cdots+\eta(c_n)\] equal
to zero.}\medskip} 

\noindent
It follows that the point \[p=\eta^{-1}(0) \in W\] is uniquely defined.
{\bf Definition:} This point \[p\]
will be called the {\bit conformal barycenter\/} of the \[c_j\]. A proof of
this lemma can easily be constructed from [Douady and Earle, \S2]. \endproof
\smallskip

{\bf Proof of 3.3.} Evidently \[h\] can be any M\"obius automorphism of \[D\]
which
carries zero to the conformal barycenter of the critical points of \[\phi\],
and which carries \[+1\] to one of the \[d\] points of \[\phi^{-1}(1)\].
\smallskip

To determine the topology of \[B'\], we proceed as follows.
Note first that the subspace \[\s_n(\bar D)\subset \s_n(\C)\],
consisting of \[n$-tuples which belong to the closed disk \[\bar D\], forms a
closed \[2n$-cell with interior equal to \[\s_n(D)\]. In fact, for each
non-zero \[\{a_1\,,\,\ldots\,,\,a_n\}\in\s_n(\C)\] consider the ray consisting
of points
\[\{ta_1\,,\, \ldots\,,\,ta_n\}\] with \[t \ge 0\]. Each such ray crosses
the boundary of \[\s_n(\bar D)\] exactly once, and the image of each such ray
in the space of ordered \[n$-tuples \[(\sigma_1\,,\,\ldots\,,\,\sigma_n)\]
crosses the unit sphere exactly once. Hence, stretching by an appropriate
factor along each such ray, we obtain the required homeomorphism from \[\s_n
 (\bar D)\] to the unit disk in \[\C^n\]. Using this construction we see also
that the subspace of \[\s_n(D)\] consisting of unordered \[n$-tuples with sum
zero is a topological \[(2n-2)$-cell.
Thus the set of Blaschke products of the form
\[\beta(z)=\mu_{a_1}(z) \cdots \mu_{a_d}(z)\] with \[a_1+\cdots+a_d=0\]
is an open topological \[(2d-2)$-cell. It follows that the collection
\[B'\] of boundary-rooted, critically
centered Blaschke products \[\beta'\] of degree \[d\] is also an open
topological \[(2d-2)$-cell.
For given any such \[\beta'\], by Lemma 3.4 there is a
unique boundary-rooted
 M\"obius automorphism \[\eta\] so that \[\beta=\beta' \circ
 \eta\] has the sum of its zeros equal to zero,
and similarly given \[\beta\] there is
a unique \[\eta\] so that \[\beta'=\beta\circ \eta^{-1}\] is critically
centered.\endproof\medskip

We will use such Blaschke products to model the map \[f\crit :
 W\crit \to W\crit\] studied in 2.3. Since \[W\crit\] is isomorphic to
a disjoint union of open disks, let us first study proper maps from such a
disjoint union to itself.\smallskip

{\bf Definition.} To any reduced mapping schema \[S=(|S|\,,\,F\,,\,w)\]
we associate the {\bit model space\/} \[B(S)\] consisting of all proper
holomorphic maps
$$ \beta:|S|\times D\;\to\; |S|\times D $$ such that \[\beta\] carries
each \[v\times D\] onto \[F(v)\times D\] by a boundary-rooted map of
degree\break \[d(v)=w(v)+1\;\] which satisfies \[\;\beta(v,0)
=(F(v)\,,\,0)\;\] whenever \[v\] is periodic under \[F\], and is critically
centered whenever \[v\] is not periodic under \[F\].\smallskip

{\QP{\bf Lemma 3.5.} \it If the schema
\[S\] has total weight \[w\], then the model space
\[B(S)\] is homeomorphic to an open cell of dimension \[2w\].\medskip}

{\bf Proof.} This follows easily from 3.2 and 3.3.\endproof

We next show that the various maps in \[B(S)\] serve as models for all possible
dynamics which can occur within the basins of attracting cycles.

{\QP{\bf Lemma 3.6.} \it Let \[M\] be a disjoint union
of finitely many open disks, and let \[\phi:M\to M\] be any proper
holomorphic map, of degree \[\ge 2\]
on each component of \[M\], such that every orbit under \[\phi\] converges
to an attracting cycle. Then \[\phi\]
is holomorphically conjugate to a map which belongs to some model space
\[ B(S)\]. Here the schema \[S\] is uniquely determined up to isomorphism;
and the number of distinct conformal isomorphisms from \[M\]
onto \[|S|\times D\] which conjugate \[\phi\] to some element of \[ B(S)\]
is equal to the order of the group \[G(S)\].\medskip}

{\bf Proof.} We may identify the complex manifold \[M\] with \[\Sigma
\times D\], where \[\Sigma\] is a finite index set. Note that \[\phi\]
extends continuously over \[\Sigma\times\partial D\]. First consider
some component \[\sigma\times D\] which is mapped to itself by some iterate
\[\phi^{\circ k}\]. Let \[d_1\cdots d_k\] be the degree
of \[\phi^{\circ k}\] on this
component. Then \[\phi^{\circ k}\] has a unique fixed point inside
\[\sigma\times D\]. After conjugating by a M\"obius automorphism
of \[\sigma\times D\], we may take this
interior fixed point to be the center point \[(\sigma\,,\,0)\]. Similarly,
there are \[\;d_1\cdots d_k-1\;\] fixed points on \[\sigma\times \partial D\],
and we can rotate \[\sigma\times \bar D\]
so that one of these fixed points is \[(\sigma\,,\,1)\].
Pushing forward, we find corresponding interior and boundary points for
each of the \[k\] disks in the cycle. Finally we work outwards, first choosing
corresponding preferred points for each of the disks which maps immediately
to a disk on this cycle, and then continuing inductively. Details
will be left to the reader.\endproof

We can sharpen 3.6 as follows.

{\QP{\bf Lemma 3.7.} \it The effective symmetry group \[\bar G(S)\] of \S2
operates smoothly on the model space \[B(S)\] in such a way
that two maps in \[B(S)\] are conformally conjugate to each other
if and only if they belong to the same orbit under this action. In fact
this action of \[\bar G(S)\] on \[B(S)\] is covered by an action of the full
symmetry group \[G(S)\] on the product \[B(S)\times|S|\times  D\] with
the following property: Each \[g\in G(S)\] carries \[(\beta\,,\,v
\,,\,z)\] to a triple of the form \[(\beta'\,,\,h(v,z))\]
where \[h=h_{g,\beta}\] is
a conformal automorphism of \[|S|\times D\] depending on \[g\] and \[\beta\],
and where \[\beta'=h\circ\beta\circ h^{-1}\].
\medskip}

{\bf Remark.} We will see in \S4 that \[B(S)\] is canonically diffeomorphic
to the model hyperbolic component \[\hl_0^S\]. Since \[\bar G(S)\] acts
linearly on \[\hl_0^S\], it is hardly surprising that it operates smoothly
on \[B(S)\].
However, to avoid a circular argument, we must prove this fact from scratch.

{\bf Proof of 3.7.} Since the action of the subgroup \[\Aut(S)\subset G(S)\]
on \[B(S)\] and on \[B(S)\times|S|\times D\] is quite easy to describe, let
us concentrate on the complementary subgroup \[N(S)\subset G(S)\]. In other
words, we will only discuss those conformal automorphisms of \[|S|\times D\]
which carry each component onto itself. If such an automorphism \[h\]
conjugates some map \[\beta\in B(S)\] into another map \[\beta'=h^{-1}\circ
\beta\circ h\] which also belongs to \[B(S)\],
then evidently \[h\] must preserve center points, and hence must
have the form\break \[h(v\,,\,z)=(v\,,\,\eta_v\,z)\]. (Here \[h\] is
the inverse of the map considered in the previous paragraph.)
Now the condition that
\[\beta'\] is boundary-rooted takes the form
$$	\beta(v\,,\,\eta_v)\;=\;(F(v)\,,\,\eta_{F(v)})\;. \eqno (1)$$
%
%
%
Evidently this condition depends on the particular map \[\beta\]
we have chosen. However, it is not difficult to check that the possible
solutions vary continuously as we vary \[\beta\]. Since the space \[B(S)\]
is contractible, this means that if we find a solution for one point \[\beta_0
\in B(S)\], then we obtain corresponding solutions for all points \[\beta\in
B(S)\].
Define the ``center point'' \[\beta_0\in B(S)\] to be the mapping
defined by 
$$ \beta_0(v,z)\;=\;f_0(v,z)\;=\;(F(v)\,,\, z^{d(v)})\;. \eqno (2)$$ 
(Compare 2.4 and 3.8.) Then it is true, almost
by definition, that the group consisting of all holomorphic conjugacies
from \[\beta_0\] to itself can be identified with the symmetry group
\[G(S)\]. Thus, deforming the solutions of (1) continuously,
we obtain an action of the group \[N(S)\], and hence of \[G(S)\], on the
space \[B(S)\]
and on \[B(S)\times|S|\times D\].\smallskip

Next we must ask which elements of the group \[N(S)\subset G(S)\] act trivially
on \[B(S)\]. Note that the automorphism \[h(v,z)=(v\,,\,\eta_v\,z)\] of \[|S|
\times D\] commutes with the map \[\beta(v,z)=(F(v)\,,\,\beta_v(z))\] if and
only if
$$	\beta_v(\eta_v\,z)\;=\;\eta_{F(v)}\beta_v(z)\;. \eqno (3)$$
In most cases, we can choose each map \[\beta_v\] so that its set of
zeros admits no non-trivial rotation. Whenever this is the case, the equation
(3) admits only the trivial solution \[\eta_v=\eta_{F(v)}=1\]. However, in
the special case where the vertex \[v\] has weight\break \[w(v)=1\], and is not
periodic under \[F\], the map \[\beta_v\] must be critically centered of
degree two, and hence must have the form
$$	\beta_v(z)\;=\;\mu_a(z)\mu_{-a}(z)\;=\;{z^2-a^2\over 1-\bar a^2z^2}
\;.$$
In this special case, the set of zeros of \[\beta_v\] necessarily
admits an \[180^\circ\] rotation. Hence equation (3) admits the solutions
\[\eta_v=\pm 1\,,\;\eta_{F(v)}=1\].
In this way, it is not difficult to check
that the subgroup \[N_0(S)\subset G(S)\] of 2.9,
which acts trivially on \[\p^S\], is precisely equal to the subgroup which
acts trivially on \[B(S)\]. Thus we obtain
the required effective action of the quotient
group \[\bar G=G/N_0\] on \[B(S)\]. The details of this argument
are not difficult, and will be left to the reader.\endproof\smallskip

Finally, we will need the following result, which is essentially well known.

{\QP{\bf Lemma 3.8.} \it Each model space \[B(S)\] contains one and only one
map \[\beta_0\] which is
post-critically finite. It is given by formula \[(2)\] above,
and is characterized by the
properties that each component \[v\times D\] contains
exactly one critical point, and that
\[\beta_0\] maps critical points to critical points.\smallskip}

{\bf Proof.} First consider a proper holomorphic map \[\beta : D
 \to D\], which has a attracting fixed point at zero. If the multiplier
\[\lambda
 =\beta'(0)\] is non-zero, then by the K\oe nigs linearization theorem
we can choose a local coordinate \[\zeta\]
in a neighborhood of zero so that \[\beta\] corresponds to the map \[\zeta
 \mapsto \lambda \zeta\].       
Let us extend this coordinate system to a conformal isomorphism
between an open set \[U \subset D\] and an open disk \[|\zeta|<r\],
with \[r\] as large as possible. Note that \[r\] cannot be infinite, since
no open subset of \[D\] is conformally isomorphic to the whole complex line.
If \[r\] is maximal, then there must
be some obstruction to extending further, and this obstruction can only be
a critical point lying on the boundary of \[U\]. Thus, under this hypothesis,
there must be a critical point in \[D\] whose forward orbit is {\it not}
eventually periodic.
 
Therefore, in the post-critically finite case, the multiplier \[\lambda\] must
be zero. Suppose then that \[\lambda=0\]. By B\"ottcher's
theorem, we can choose the
local coordinate \[\zeta\] so that \[\beta\] corresponds to the map \[\zeta
 \mapsto \zeta^k\] for some \[k \ge 2\].        
Again we extend to a conformal isomorphism from \[U \subset D\] onto the open
disk \[|\zeta|<r\] with \[r\] as large as possible. If the maximal \[r\]
satisfies \[r<1\], then again there must be a critical point on the boundary
of \[U\], whose forward orbit is not eventually periodic. Thus, if all
critical orbits are eventually periodic, it follows that \[r=1\]. But then
the boundary of \[U\] maps into itself under \[\beta\], and it follows
easily that \[U=D\],
so that \[\beta\] is conformally conjugate to the map \[z \mapsto z^k\]. In
fact if \[\beta(1)=1\] then \[\beta(z)=z^k\]. In this case, note that the
unique critical point of \[\beta\] is the {\it only\/} point of \[D\] whose
forward orbit is eventually periodic.

To complete the proof, we must also consider the case of a critically
centered Blaschke product \[\beta : D \to D\], with \[\beta(1)=1\]. If all
of the critical points of \[\beta\] map to zero, then we must show again that
\[\beta(z)=z^d\]. But this is clear, since \[D\] is exhibited as a branched
covering of \[D\] with unique branch point at the critical value zero.
\endproof
\bigskip\bigskip

\centerline{\bf \S4. Hyperbolic Components are Topological Cells.}
\bigskip

The object of this section is to prove the following two results.

{\QP{\bf Theorem 4.1.} \it If \[H_\alpha\subset\cl^{S_0}\] is any hyperbolic
component whose elements \[f\] have reduced mapping schema \[\bar S(f)\]
isomorphic to \[S\], then \[H_\alpha\] is diffeomorphic to the model space
\[B(S)\]. This diffeomorphism is canonically defined, up to composition
with an element of the group \[\bar G(S)\] which acts on \[B(S)\].\medskip}

\noindent The proof will also demonstrate the following.

{\QP{\bf Corollary 4.2 (McMullen).} \it Each hyperbolic component \[H_\alpha\] contains
one and only one map \[f_\alpha\] which is post-critically finite.\medskip}

Equivalently, such a ``center'' map \[f_\alpha\]
has the property that each component of the complement of \[J(f_\alpha)\]
contains one and only one pre-critical point. \medskip


To begin the argument, suppose that \[f\] is a hyperbolic map
belonging to some connectedness
locus \[\cl^{S_0}\], with reduced mapping schema \[\bar S(f)\]
isomorphic to \[S\]. Then the open set \[W\crit\subset K(f)\], as defined
in 2.3, satisfies the hypothesis of Lemma 3.6. That is, \[W\crit\] is a union
of finitely many components, each conformally isomorphic to \[D\], and the
first return map \[f\crit\] from \[W\crit\] to itself has degree at least
two on each component. Hence, by 3.6, there exists a conformal isomorphism
\[h:W\crit\to |S|\times D\] which conjugates \[f\crit\]
to some map \[\beta\in B(S)\].\smallskip

However, we must
be careful since \[h\] and \[\beta\] are
not uniquely defined. (Compare 3.7.) To
deal with this non-uniqueness, we proceed as follows. Note first that for
each component\break \[W_\alpha\subset W\crit\] the closure \[\bar W_\alpha\]
is homeomorphic to the closed unit disk \[\bar D\], in a homeomorphism which
is conformal throughout the interior.  (See [DH1, pp. 13, 24, 26].)\smallskip

{\bf Caution.} The closed disks \[\bar W_v\] are not
always disjoint from each other, so the closure of \[W\crit\] need
not be homeomorphic to \[|S|\times \bar D\]. 
 Since boundary points of
\[W\crit\] will play an important role in our discussion, it will be
convenient to introduce the following notation. {\it Let \[\hat W\crit=\coprod
 \bar W_v\] be the disjoint union of the closures of the components of
\[W\crit\], and let \[\hat f\crit:\hat W\crit\to\hat W\crit\] be the continuous extension
of \[f\crit\].}

In order to choose some specific conformal isomorphism between \[|S|\times
D\] and \[ W\crit\], we must first choose some isomorphism \[\iota:S\to
\bar S(f)\]. Evidently the number of ways of doing this
is equal to the order of the automorphism group \[{\Aut}(S)\]. For each
component \[ W_{\iota(v)}\] of \[ W\crit\], we must then
choose one boundary point \[q(v)\]
which will correspond to the boundary base point
\[(v,1)\in v\times \bar D\].
These boundary points are to be chosen as follows.\smallskip

{\bf Definition 4.3. Boundary markings.} Let \[f\in \hl^{S_0}\] be a hyperbolic
map with reduced schema \[\bar S(f)\] isomorphic to \[S\]. By a {\bit
 boundary marking\/} \[q:|S|\to\partial\hat W\crit\] of \[f\], {\bit covering
the isomorphism\/} \[\iota:S\to\bar S(f)\],
we will mean a function which assigns to each vertex \[v\in |S|\] a boundary
point \[q(v)\in\partial W_{\iota(v)}\],
satisfying the condition that\break \[\hat f\crit(q(v))=q(F(v))\].

{\QP {\bf Lemma 4.4.} \it Such boundary markings always exist. In fact,
the number of distinct boundary markings for \[f\] is precisely equal to the
order of the symmetry group \[G(S)\]. Given such a boundary marking \[q\],
there
is one and only one homeomorphism \[\hat q:|S|\times \bar D\buildrel\cong\over
\longrightarrow \hat W\crit\], conformal throughout the interior, which
satisfies \[\hat q(v,1)=q(v)\], for which the map \[\;\beta\;=\;\hat q^{-1}\circ
f\crit\circ \hat q\,:\,|S|\times D\,\to\,|S|\times D\;\] belongs to the
model space \[B(S)\].\medskip}

The proof is completely analogous to the proof of 3.6, and will be omitted.
\endproof

 Note also that as we deform the hyperbolic map \[f\] within a
small neighborhood, there is a corresponding deformation of any given
boundary marking. This follows easily from the theorem of Ma\~n\'e-Sad-Sullivan
and Lyubich, which asserts that the entire Julia set \[J(f)\] varies
continuously as we vary the hyperbolic map \[f\].
We can now make a more precise restatement of 4.1:
\eject

{\QP{\bf Theorem 4.1$\bf '$.} \it If \[H_\alpha\subset \cl^{S_0}\]
is any hyperbolic
component whose elements\break \[f\in H_\alpha\] have mapping
schema \[\bar S(f)\]
isomorphic to \[S\], then:

\ref$(1)$ we can choose a boundary marking \[\;\;
q_f:|S|\to\partial\hat W\crit(f)\;\;\]
which varies\break continuously as \[f\] varies over \[H_\alpha\],

\ref$(2)$ there is an associated extension
\[\;\;	\hat q_f:|S|\times\bar D \buildrel \cong\over\longrightarrow
\hat W\crit\; \;\]
as in 4.4, and 

\ref$(3)$ the resulting
correspondence \[\;~f\,\mapsto \,\hat q_f^{-1}\circ\hat f\crit\circ \hat q_f\,
\in\, B(S)~\;\] maps the comp\-onent \[H_\alpha\] diffeomorphically onto the
space \[B(S)\] of Blaschke products.\bigskip}

The proof is a generalization of that given by Douady and Hubbard in the
quad\-ratic case. However, more care is needed, since it is necessary to
make a choice of boundary markings. For example, a priori some component
\[H_\alpha\] might have a non-trivial fundamental group. If this were the case,
then starting with a boundary marking \[q_0\] for \[f_0\] and deforming it
as we follow a loop around \[H_\alpha\] we might end up with a different
boundary marking for \[f_0\].
In order to prove that this cannot happen, We proceed as follows.\smallskip

Let \[\tilde H_\alpha\] be the space consisting of all pairs
\[(f,\,q)\] where \[f\] belongs to \[H_\alpha\] and \[q\] is a boundary
marking for \[f\]. This space \[\tilde H_\alpha\] has a natural topology,
and there is a natural projection \[(f,\,q)\mapsto f\] from \[\tilde
 H_\alpha\] to \[H_\alpha\].
Define a map \[\Phi:\tilde H_\alpha \to B(S)\;\] by the construction
$$ \Phi:(f,\,q)~ \mapsto ~\hat q^{-1} \circ \hat f\crit \circ \hat q~~, $$
with \[\hat q\] as in 4.4. We will prove the following.

{\QP{\bf Lemma 4.5.} \it The space \[B(S)\] is evenly
covered under this map\vskip -.1in
 $$\Phi:\tilde H_\alpha \to B(S)\;.$$\par}

 Since \[B(S)\] is simply connected by 3.5, this implies
that there is a section\break \[B(S)\to\tilde H_\alpha\].
The composition \[B(S)\to\tilde H_\alpha\to
H_\alpha\] is then the required diffeomorphism from \[B(S)\]
onto \[H_\alpha\], with inverse as described in 4.1$'$.\smallskip

{\bf Proof of 4.5.}
Let us start with a map \[f_0\in H_\alpha\] with
boundary marking \[q_0\] and with associated model map \[\;b_0=\hat q_0^{-1}
 \circ \hat f_0\crit \circ \hat q_0\in B(S)\,\]. Choose two radii \[r_1<r_2<1\]
so that

(1) every critical point of the map \[b_0:|S|\times D\to|S|\times D\]
is contained in the union \[|S|\times D_{r_1}\] of disks of radius \[r_1\],
and

(2) so that \[b_0\] maps \[|S|\times\bar D_{r_2}\] into
\[|S|\times D_{r_1}\].

\noindent Let \[U\subset B(S)\] be a simply connected neighborhood
of \[b_0\] which is small enough so that all
maps \[b\in B(S)\] will satisfy these same conditions. That is, the union
\[|S|\times D_{r_1}\] must contain all critical points of \[b\],
and must contain
the image \[b(|S|\times\bar D_{r_2})\]. Then we will construct a new polynomial
map \[f_b\in H_\alpha\] by quasi-conformal surgery.\smallskip

First construct a new function \[b'\] from \[|S|\times D\]
to itself as follows.
Let \[b'\] coincide with \[b\] on \[|S|\times \bar D_{r_1}\] and with \[b_0\]
outside of \[|S|\times D_{r_2}\]. Now interpolate linearly
on each intermediate region \[r_1\le |z|\le r_2\], setting
$$	b'(v,z)\;=\;t\,b_0(v,z)+(1-t)b(v,z)\;, $$
where \[t=(|z|-r_1)/(r_2-r_1)\]. We will assume that
the neighborhood \[U\] is small enough so that this linear interpolation
yields a new function \[b'\] which is a quasi-conformal local homeomorphism
throughout these annuli \[r_1\le |z|\le r_2\]. Note that \[b'\] is actually
holomorphic outside of these annuli.\smallskip

We now follow the surgery procedure originated by Douady and Hubbard.
(Compare [D1], [DH1], [McM], [Sh1], [Sh2], [Su3].)
Identifying \[|S|\times D\] with \[W\crit\]
under \[\hat q_0\], we obtain a quasi-conformal map \[g_b\] from the Riemann
sphere to itself by setting \[g_b(z)=f(z)\] outside of \[W\crit\], with
\[g_b=\hat q_0\circ b'\circ\hat q_0^{-1}\] within \[W\crit\]. Then \[g_b\]
is holomorphic except on a collection of annuli, one of which lies in each
component of \[W\crit\]. Since every orbit under \[g_b\] passes through these
bad annuli at most once, we can pull back the standard conformal structure
\[\mu^0\] under the iterates of \[g_b\] to obtain a new conformal structure
\[\mu_b\] on \[\hat\C\] which is invariant under \[g_b\]. Using the measurable
Riemann mapping theorem, we see that there is a quasi-conformal automorphism
\[h_b\] of the Riemann sphere which transforms \[\mu^0\] to \[\mu_b\].
(Compare Ahlfors \& Bers or Lehto \& Virtanen.) Thus the conjugate mapping
\[f_b=h_b^{-1}\circ g_b\circ h_b\] preserves the standard structure \[\mu^0\],
or in other words is holomorphic. Furthermore, if we choose \[h_b\] to be
doubly tangent to the identity at infinity, then it is unique, and
depends (real) differentiably on the parameter \[b\]. Hence the holomorphic
map \[f_b\] also depends differentiably on \[b\]. It is not difficult to check
that \[f_b\] is a monic centered polynomial with a preferred boundary marking.
Thus it belongs to our space \[\hat H_\alpha\], and
we have constructed a smooth map \[b\mapsto f_b=s(b)\] from the open
set\break \[U\subset B(S)\] to \[\hat H_\alpha\].\smallskip

We would like to prove that \[s\] is a local section of the projection map\break
\[\Phi:\hat H_\alpha\to B(S)\]. That is, we would like to prove
that the composition $$\Phi(s(b))\;=\;\Phi(f_b)\,\in\, B(S)$$ is equal to \[b\].
Examining the construction, we certainly see that if we restrict \[b\] to
the disks \[|S|\times D_{r_1}\] and restrict \[\Phi(f_b)\] to a corresponding
neighborhood of its critical set, then the two are holomorphically conjugate.
It is not too difficult to show that
this holomorphic conjugacy extends inductively over the iterated pre-images
of \[|S|\times D_{r_1}\], since each of these can be considered as a branched
covering of the previous one. Passing to the direct limit, we see that \[b\]
is holomorphically conjugate to \[\Phi(f_b)\] on the entire space
\[|S|\times D\].\smallskip

Thus we have constructed a smooth map
\[b\mapsto s(b)=f_b\] from an open set\break \[U\subset B(S)\] to \[\tilde H_\alpha\]
such that the image \[\Phi\circ s(b)=
\Phi(f_b)\in B(S)\] is holomorphically conjugate to
\[b\]. Hence \[\Phi\circ s(b)=g_b( b)\] for some group element \[g_b\in\bar G(S)\].
It follows easily from smoothness that \[g_b\] can be chosen as a
constant, independent of \[b\]. Then \[s'=s\circ g_b^{-1}\] is the required
local section of the projection \[\Phi\]. This completes the proof of 4.5,
$4.1'$ and 4.1.\endproof\medskip

{\bf Proof of 4.2.} This follows from the argument above, together with 3.8.
\endproof

\end

\def\endproof{  \rlap{$\sqcup$}$\sqcap$ \smallskip}
\def\C{{\bf C}}
\def\QP{\medskip\leftskip=.4in\rightskip=.4in\noindent}
\def\[{$\,}
\def\]{\,$}
\font\bit=cmssi12 at 12truept
\def\Aut{ {\rm Aut}}
\def\crit{^C}
\def\cl{{\cal C}}
\def\hl{{\cal H}}
\def\p{{\cal P}}
\font\tenmsy=cmbsy10
\textfont8=\tenmsy
\mathchardef\ssm="7872
\pageno=19

\centerline{\bf \S5. Analytic isomorphism between hyperbolic components.}
\medskip

If \[H_\alpha \subset \cl^{S_1}\] and \[H_\beta \subset \cl^{S_2}\] are
two different hyperbolic components with reduced mapping schema
isomorphic to \[S\], then by Theorem 4.1 there are diffeomorphisms
$$	H_\alpha \to B(S) \leftarrow H_\beta \,,	$$
uniquely defined up to a choice among finitely many boundary markings, or
equivalently up to the action of the group \[\bar G(S)\] on \[B(S)\]. The
composition mapping \[H_\alpha\] to \[H_\beta\] will be called the
{\bit canonical diffeomorphism\/} between these two sets. We
will prove the following.

{\QP {\bf Theorem 5.1.} {\it This canonical diffeomorphism \[H_\alpha \to
 H_\beta\] between open subsets of complex affine spaces is biholomorphic.}
\smallskip}

That is, it is holomorphic, with holomorphic inverse.
In particular, it follows that the canonical diffeomorphism from \[H_\alpha\]
to the standard model \[H_0^S\] is biholomorphic. Note that this
diffeomorphism is unique up to the action of the finite group
\[\bar G(S)\] of linear automorphisms of \[H_0^S\]. The proof will
be based on the following.\smallskip

{\bf Definition.} We will say that a map \[f \in \p^{S_1}\] satisfies
a {\bit critical orbit relation\/} if either (1) the \[w\] critical points
of \[f\] are not all distinct, or (2) the associated critical orbits are not
disjoint from each other, or (3) some critical orbit is eventually
periodic. It is not difficult to show that
the set of all \[f\] which satisfy some critical orbit
relation forms a countable union of algebraic varieties in the affine space
\[\p^{S_1}\]. 

{\QP {\bf Lemma 5.2.} {\it If \[f_1 \in H_\alpha\] does not satisfy any
critical orbit relation, then
the\break canonical diffeomorphism from \[H_\alpha\] to \[H_\beta\] is biholomorphic
throughout a neighborhood of \[f_1\].}	\smallskip}

{\bf Proof.} First consider the mapping schema \[S=(w)\], with a single vertex
of degree \[d=w+1\]. Then every
\[f_1\in H_\alpha\] has a unique attracting periodic orbit. To simplify
the notation, let us assume that it is a fixed point \[p_1=p(f_1)\]. We can
construct a local holomorphic coordinate system for \[H_\alpha\] near \[f_1\]
as follows. Since there are no critical orbit relations,
the multiplier \[\lambda_1=f'_1(p_1)\in D\] cannot be zero. Hence every
\[f\] close to \[f_1\] has a fixed point \[p=p(f)\] close to \[p_1\] with
multiplier \[\lambda=\lambda(f) \in D\ssm \{0\}\]. Let us choose a K{\oe}nigs
linearizing
coordinate \[z \mapsto \zeta(z)\] near \[p(f)\] so that \[\zeta(f(z))=\lambda
 \zeta(z)\]. Evidently \[\zeta\] extends to a holomorphic
map from the basin of attraction \[W(f)\] into \[\C\] satisfying this same
identity \[\zeta(f(z))=\lambda
 \zeta(z)\]. Since there are no critical orbit relations,
the values \[\zeta(c_1)\,,\,\ldots\,,\,\zeta(c_w)\] at
the critical points of \[f\] are all distinct and non-zero. In fact no
ratio \[\zeta(c_j)/\zeta(c_k)\] can be equal to a power of \[\lambda\].
Let us normalize so that \[\zeta(c_1)=1\]. {\it Then the numbers \[\lambda\,
 ,\, \zeta(c_2)\,,\,\ldots\,,\, \zeta(c_w)\] form the required local
holomorphic coordinates near \[f_1\].}\smallskip

To see that the correspondence \[f \mapsto \zeta(c_j)\] is holomorphic,
note that \[c_j=c_j(f)\] depends holomorphically on \[f\], and hence that
$$	\zeta(c_j)= {\lim \limits_ {k \to \infty} }(f^{\circ k}(c_j)-p)/(
	f^{\circ k}(c_1)-p)	$$
is a uniform limit of holomorphic functions. We must show that these
coordinates,\break in a sufficiently small neighbood of \[f_1\], determine \[f\]
uniquely. Choose open sets\break \[U_0 \subset U_1 \subset \cdots \subset
 \hat W(f)\] as follows. Choose some number \[0<r<1\] which is smaller than all
of the \[|\zeta(c_j)|\], and let \[U_k\] be the component containing
the fixed point \[p\] in the region \[|\lambda^k\zeta(z)|<r\]. Then each
\[U_{k+1}\] can be described as a branched covering of \[U_k\], with the
restriction of \[f\] to \[U_{k+1}\] as projection map,
and with the critical points of \[f\] in \[U_{k+1}\] as branch points.
Note that \[U_0\] can be identified with the open
disk \[\{\zeta:|\zeta|<r\}\]. This entire sequence
of Riemann surfaces, together with the holomorphic functions \[f\] and
\[\zeta\] on the \[U_k\], and the embedding of \[U_k\] into \[U_{k+1}\]
can be constructed inductively, starting with \[U_0\], if we are
given the numbers
\[\lambda\,,\, \zeta(c_2)\,,\,\ldots\,,\, \zeta(c_w)\] and \[r\] together
with the topological data which specifies which covering to take. Evidently
this topological data (eg. certain subgroups of finite index in
fundamental groups of punctured disks)
varies continuously with \[f\]. Passing to the direct limit of \[U_k\] as
\[k \to \infty\], we have constructed  a Riemann surface conformally
isomorphic to \[W\crit(f)\cong D\], together with the
map \[f\crit\] on \[W\crit(f)\], from the given data. By Theorems 3.6 and
4.1, this
determines \[f\] uniquely, up to the action of the effective symmetry group.

The proof for an arbitrary connected mapping schema \[S\] is similar. If
the unique attracting periodic orbit has period \[m\], then it is convenient
to write the multiplier as \[\lambda^m\], and to choose \[\zeta:W(f) \to
 \C\] so that \[\zeta(f(v',z) =\lambda\zeta(v',z)\]. Proceeding as above,
from the data \[\lambda\,,\,\zeta(c_2)\,,\,\ldots\,,\,\zeta(c_w)\,,\, r\]
and \[S\] we
can first build up a copy of the component \[W_v\] containing one of the
periodic points, and then build up the rest of \[W\crit\] inductively.
The case of an \[S\] with several connected components now follows easily.
Since these same local holomorphic coordinates can be used either in \[H_
 \alpha\] or in \[H_\beta\], this proves Lemma 5.2.\endproof

It now follows easily that this diffeomorphism \[H_\alpha \to H_\beta\] is
holomorphic everywhere. In fact the Cauchy-Riemann equations, which are
necessary and sufficient conditions for a $C^1$-smooth map to be holomorphic,
are satisfied except on a countable union of proper algebraic subvarieties,
intersected
with \[H_\alpha\]. Since such a countable union is nowhere dense in
\[H_\alpha\], Theorem 5.1 follows by continuity.\endproof


\vfil\eject
\end


\input fig.tex
\def\endproof{  \rlap{$\sqcup$}$\sqcap$ \smallskip}
\def\C{{\bf C}}
\def\R{{\bf R}}
\def\QP{\medskip\leftskip=.4in\rightskip=.4in\noindent}
\def\[{$\,}
\def\]{\,$}
\def\Aut{ {\rm Aut}}
\def\crit{^{ C}}
\def\cl{{\cal C}}
\def\hl{{\cal H}}
\def\p{{\cal P}}
\def\mixarr{\;\hbox{\raise 2pt \hbox{$\leftarrow$}}
        \hbox{\lower 1pt \hbox{\hskip -9pt $\to$}}\;}
\font\bit=cmssi12 at 12truept
\pageno=21

\centerline{\bf \S6. Real forms.}\smallskip

First consider a real polynomial map \[f_\R:\R\to\R\], of degree \[d
\ge 2\]. We can extend uniquely
to a complex polynomial map \[f:\C\to\C\]. This
\[f\] commutes with the complex
conjugation operation, which we denote by \[\gamma_0(z)=\bar z\]. If
\[f\] is hyperbolic, then as in \S2 we form the union \[ W\crit(f)\] of those
components of the attractive basin which contain critical points. Evidently
\[\gamma_0\] is an antiholomorphic mapping which carries
this set \[W\crit(f)\] onto itself. Note that \[\gamma_0\]
may permute the various components of \[W\crit(f)\;;\;\] for the critical
points of \[f\] need not be real and in fact \[W\crit(f)\] may not
contain any real points at all. In order to find an appropriate universal
model for such behavior, we consider the following construction.\smallskip

Let \[M\] be a finite union of copies of \[\C\], and let
\[\gamma'\] be any antiholomorphic involution of \[M\]. If \[f'\]
is a nowhere linear polynomial map from \[M\] to itself which commutes with
\[\gamma'\], then as in 2.7 conjugation by some affine isomorphism
\[h:M\buildrel\cong\over\longrightarrow|S|\times\C\] will carry \[f'\] to
a map  \[f=h\circ f'\circ h^{-1}\] from \[|S|\times\C\] to itself which
is monic and centered on each component \[v\times\C\]. Hence \[f\]
belongs to an appropriate space \[\p^S\] of monic centered maps.
This same affine conjugation will carry \[\gamma'\] to some
antiholomorphic involution \[\gamma=h\circ\gamma'\circ h^{-1}\]
of \[|S|\times\C\] which commutes with \[f\]. {\bf Definition:}
We will say that \[f\] belongs to the
``real subspace'' \[\p_\R^S(\gamma)\],
consisting of polynomials in \[\p^S\] which commute with \[\gamma\].
\smallskip

We can fit this construction into a group theoretic framework as
follows. Recall that the group \[G(S)\] of 2.5 can be described
as the set of all holomorphic automorphisms \[g\] of the manifold
\[|S|\times\C\] which commute with the special
map \[f_0^S(v,z)=(F(v)\,,\,z^{d(v)})\]. These automorphisms \[g\] have
the property that for any \[f \in \p^S\] the conjugate
\[g\circ f\circ g^{-1}\] will also belong to \[\p^S\]. Let us
extend to a larger group \[\hat G(S)\] by allowing also antiholomorphic
automorphisms which commute with \[f_0^S\], or equivalently which
conjugate \[\p^S\] onto itself. Then it is easy to
show that \[\hat G(S)\] is the split extension of its normal subgroup
\[G(S)\] by the two element group \[\{1,\gamma_0\}\] generated by the standard
involution \[\gamma_0(v,z)= (v, \bar z)\]. This involution
commutes with the elements of the subgroup \[\Aut(S)\subset G(S)\]. However,
for \[g\] in the normal abelian subgroup \[N(S)\subset G(S)\],
conjugation by \[\gamma_0\] carries \[g\] to \[g^{-1}\]. We can also
describe \[\hat G(S)\] as a split extension
$$ 1\to N(S)\to \hat G(S)\buildrel{\displaystyle\leftarrow}\over\rightarrow\big(
\{1,\gamma_0\} \times\Aut(S)\big)\to 1\;. $$
Proofs are easily supplied.
Thus we are led to the following.	\smallskip

{\bf Preliminary Definition 6.1.} A {\bit real form} of the mapping schema \[S\]
is an antiholomorphic involution \[\gamma:|S|\times\C\to|S|\times\C\] which
commutes with \[f_0^S\], and hence belongs to the group
\[\hat G(S)\]. Two such involutions (ie., two such real forms)
are to be considered as {\bit isomorphic\/} if they belong to the
same conjugacy class in \[\hat G(S)\]. Each real form \[\gamma\] gives
rise to an anti-linear involution \[f \mapsto \gamma\circ f\circ\gamma\]
of the complex \[w$-dimensional affine space \[\p^S\]. The fixed point set
of this involution is a real \[w$-dimensional affine space \[\p_\R^S(
 \gamma)\], consisting of maps in \[\p^S\] which commute with \[\gamma\].
We will call \[\p_\R^S(\gamma)\]
the {\bit real form} of \[\p^S\] associated with
\[\gamma\]. If we intersect this real affine space with the
connectedness locus \[\cl^S\], then we obtain a set \[\cl_\R^S(\gamma)\] which
is called a {\bit real connectedness locus\/}, or a {\bit real form\/}
of \[\cl^S\]. Similarly,
the intersection \[\hl_\R^S(\gamma)=\hl^S\cap\p^S_\R(\gamma)\] is the
associated {\bit real hyperbolic locus\/}, and its connected components are
the associated {\bit real hyperbolic components\/}. We will see in 6.4
that these real hyperbolic components are just the intersections \[H_\alpha\cap
\p_\R^S(\gamma)\], where \[H_\alpha\] can be any component of \[\hl^S\]
which intersects
\[\p_\R^S(\gamma)\], or equivalently whose center point \[f_\alpha\] commutes
with \[\gamma\]. As a special case, we can consider the {\bit principal
hyperbolic component\/} \[H_{0\R}^S(\gamma)=H_0^S\cap\p_\R^S(\gamma)\].
\smallskip

The {\bit real symmetry group} \[G_\R(S,
 \gamma)\] is defined to be the centralizer of \[\gamma\] in \[G(S)\],
that is the subgroup consisting
of all elements \[g \in G(S)\] which commute with \[\gamma\].
Evidently %
each \[g\in G(S)\] acts linearly on the spaces
\[\p_\R^S(\gamma) \supset
 \cl_\R^S(\gamma) \supset \hl_\R^S(\gamma)\supset H_{0\R}^S(\gamma)\]
by the correspondence \[f\mapsto g\circ f\circ g^{-1}\]. Hence
an appropriate quotient group \[\bar G_\R(S,\gamma)\subset \bar G(S)\]
acts effectively on these spaces.
\smallskip

Similarly, if \[M\] is a finite union of copies of the unit disk
\[D\], and \[\gamma'\]
is an antiholomorphic involution of \[M\], then we can consider proper
holomorphic maps \[f':M\to M\] which commute with \[\gamma'\], and have
degree \[\ge 2\] everywhere. As in 3.6, we can identify \[M\]
with some \[|S|\times D\] so that \[f'\] corresponds to a map \[f\in B(S)\].
Furthermore, as in 3.7, the group \[\hat G(S)\] acts on \[B(S)\times
|S|\times D\]. Thus \[\gamma'\] corresponds to an antiholomorphic involution
\[\gamma\in\hat G(S)\],
and \[f\] belongs to the space \[B_\R(S,\gamma)\] of maps which commute
with \[\gamma\]. In particular, if \[f'\] is any hyperbolic map which
belongs to some locus \[\hl_\R^{S_1}(\gamma_1)\], then the space \[M
=W\crit(f')\] is isomorphic to a finite union of disks, and is mapped into
itself by the antiholomorphic involution \[\gamma_1\], so the above discussion
applies. If the corresponding map
\[f:|S|\times D\to |S|\times D\] belongs to the model space
\[B_\R(S,\gamma)\], then we will say that the real hyperbolic
component containing \[f'\] has the {\bit type\/} \[(S,\gamma)\].\smallskip

In practice, a slightly different construction often seems more natural.
Again suppose that \[M\] is a finite union of copies of \[\C\], that \[f':
M\to M\] is a nowhere linear proper holomorphic map, and that \[\gamma'\]
is an antiholomorphic involution 
commuting with \[f'\].

{\QP{\bf Lemma 6.2.} \it We can choose a conformal isomorphism \[h:M\cong|S|
\times\C\] which conjugates \[\gamma'\] to a standard involution of the
form \[\;\gamma:(v,z)\mapsto(\bar v, \bar z)\,\], and which conjugates
\[f'\] to a map \[f:|S|\times\C\to|S|\times\C\] which is centered with
leading coefficient \[\pm 1\] on each component \[v\times\C\].\medskip}

Here the vertex map \[v\leftrightarrow\bar v\] can be any
element of order \[\le 2\] in the group \[\Aut(S)\].
Thus we can put the involution \[\gamma\] into a convenient
standard form, at the cost of allowing both \[+1\] or \[-1\] as leading
coefficients for \[f\]. The proof is not difficult.\endproof

With this formulation, the real form of \[S\] is described by the involution
\[v\leftrightarrow\bar v\], together with the leading coefficient function
\[|S|\to\{\pm 1\}\]. However, one then needs to do a little work to
decide when two such real forms are isomorphic to each other.

As an example, consider the space \[\p^{(2)}\] of monic centered cubic
polynomials. Using the formulation 6.1, we note that
\[\hat G(2)\] is the direct sum of two cyclic
groups of order two, and it follows easily that there are two suitable
antiholomorphic involutions, namely the standard involution
\[\gamma_0:
z\mapsto \bar z\] and the exotic involution \[\gamma_1:
z\mapsto -\bar z\].\break Thus we consider the two spaces \[\p_\R^2(\gamma_j)\]
consisting of monic centered cubic maps satisfying \[\overline{f( z)}=(-1)^j
f(\bar z)\], for \[j=0,\,1\].
Using the formulation 6.2, we consider instead the two spaces \[\p_\R^{2^\pm}\]
consisting of centered cubic maps with real coefficients, and with leading
coefficient \[\pm 1\]. (Compare Figure 3.)\smallskip


\pageinsert
\insertRaster 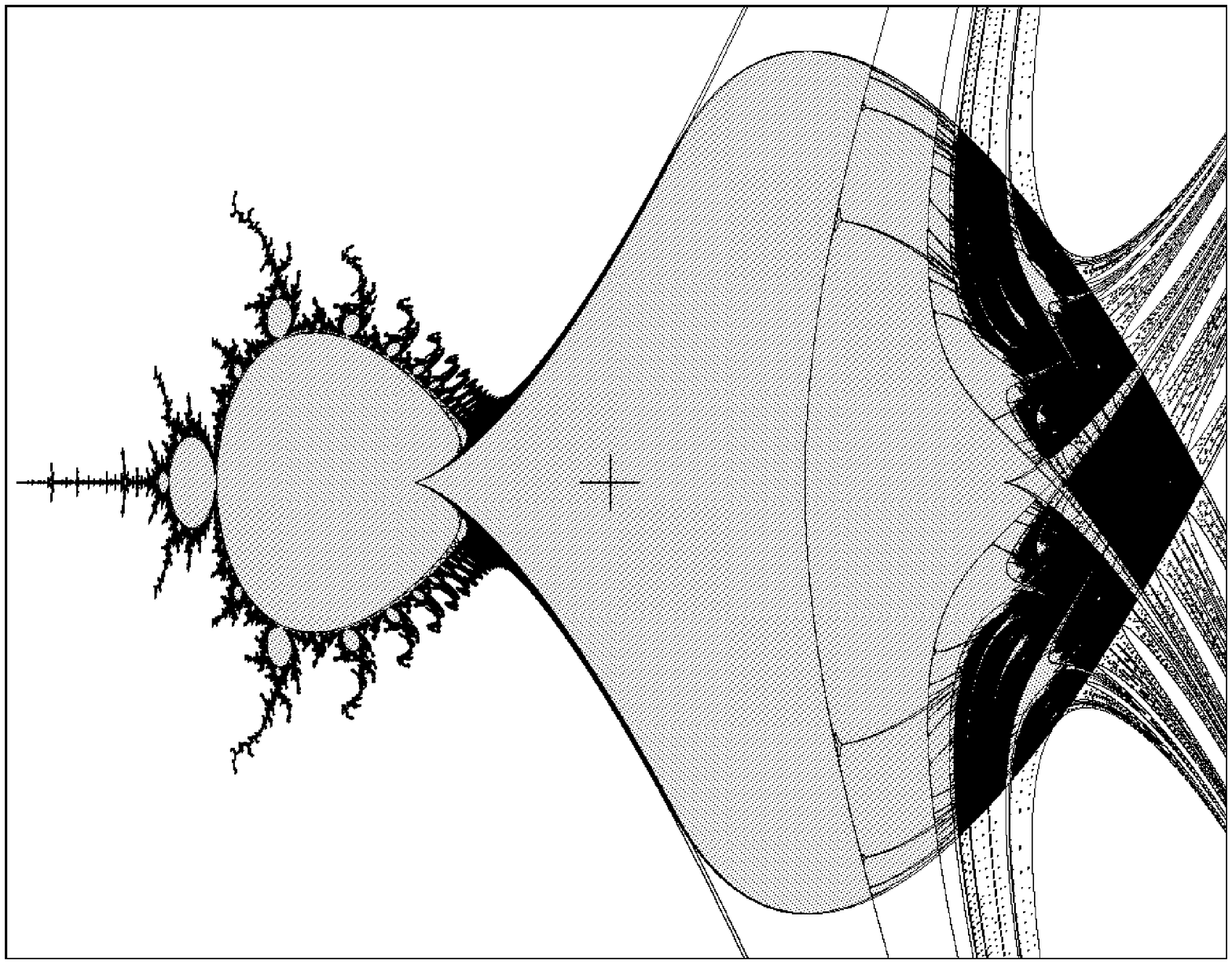 pixels 1152 by 899 scaled 225
\smallskip
\insertRaster 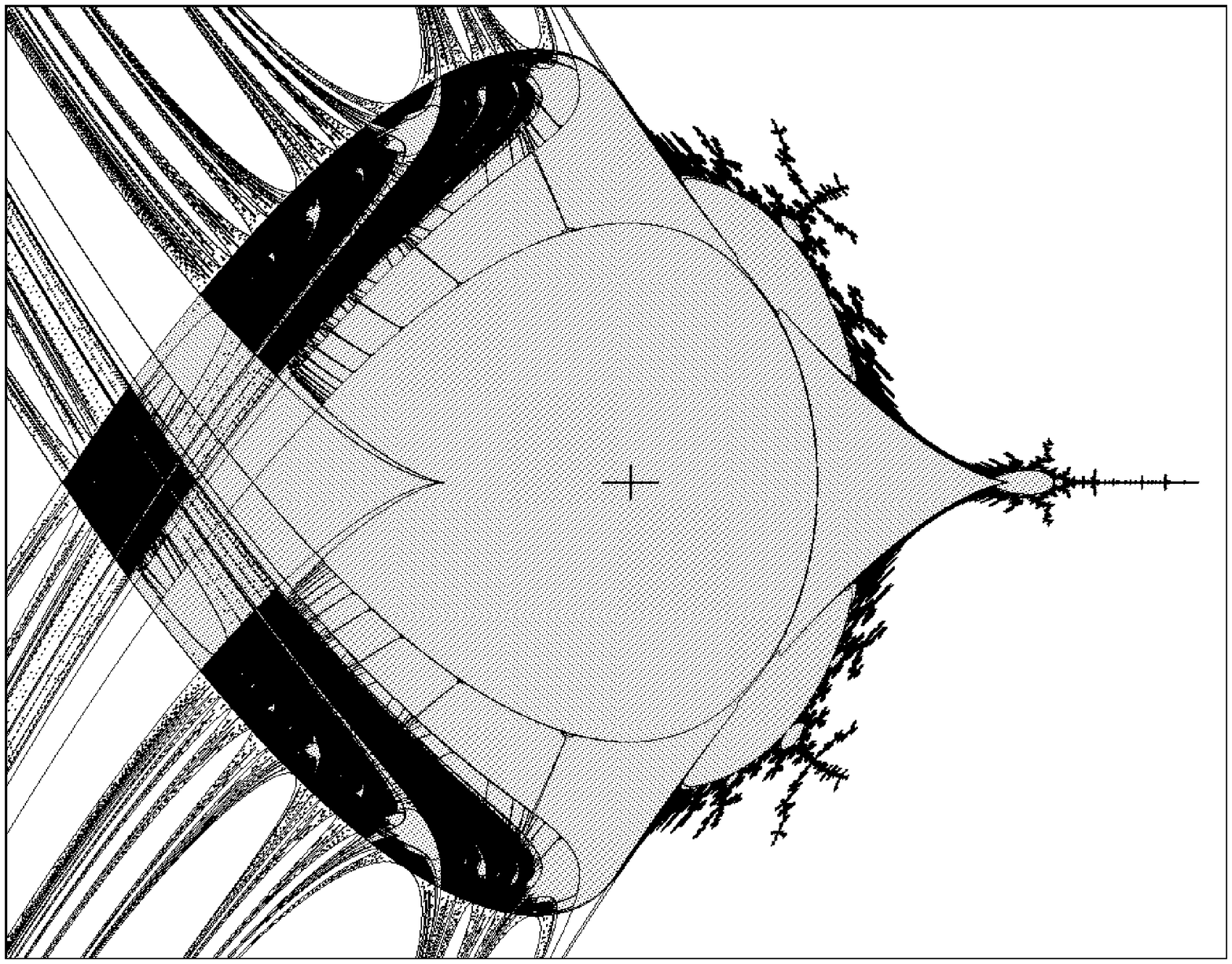 pixels 1152 by 899 scaled 225
{\QP \figf Figure 3. The spaces \[\p_\R^{(2)^+}\] and \[\p_\R^{(2)^-}\] of real
cubic maps, intersected with the complex connectedness locus \[\cl^{(2)}\].
These pictures show the \[(A,b)$-plane where \[x\mapsto \pm x^3-3Ax+b\].
(Compare [M3].)\par}
\endinsert

	\smallskip

There is a completely analogous discussion for any
even \[w\ge 2\]. The space \[\p^{(w)}\] of polynomials
of degree \[w+1\] has exactly two real forms, up to isomorphism. In fact
\[\hat G(w)\] is a dihedral group of order \[2w\] with two distinct conjugacy
classes of anti-holomorphic involutions. It is most convenient to identify the
corresponding real forms with the spaces \[\p_\R^{(w)^+}\] and \[\p_\R^{
(w)^-}\], consisting of real centered polynomials with leading coefficient
\[+1\] or \[-1\] respectively. The associated groups \[G_\R((w)^\pm)
=\bar G_\R((w)^\pm)\] are all cyclic of order two, generated by the symmetry
\[f(z)\mapsto -f(-z)\]. This symmetry
is visible as the reflection \[(A,b)\mapsto (A,-b)\] in Figure 3.
On the other hand, if \[w\] is odd, then the dihedral group \[\hat G(w)\]
has just one conjugacy class of antiholomorphic involutions. Hence there is
a unique real form, consisting of real monic centered polynomials. In this
case, the symmetry group \[G_\R((w)^+)\] is trivial.\smallskip

\midinsert
\hbox to \hsize{
\hfil
\hbox{\RasterBox {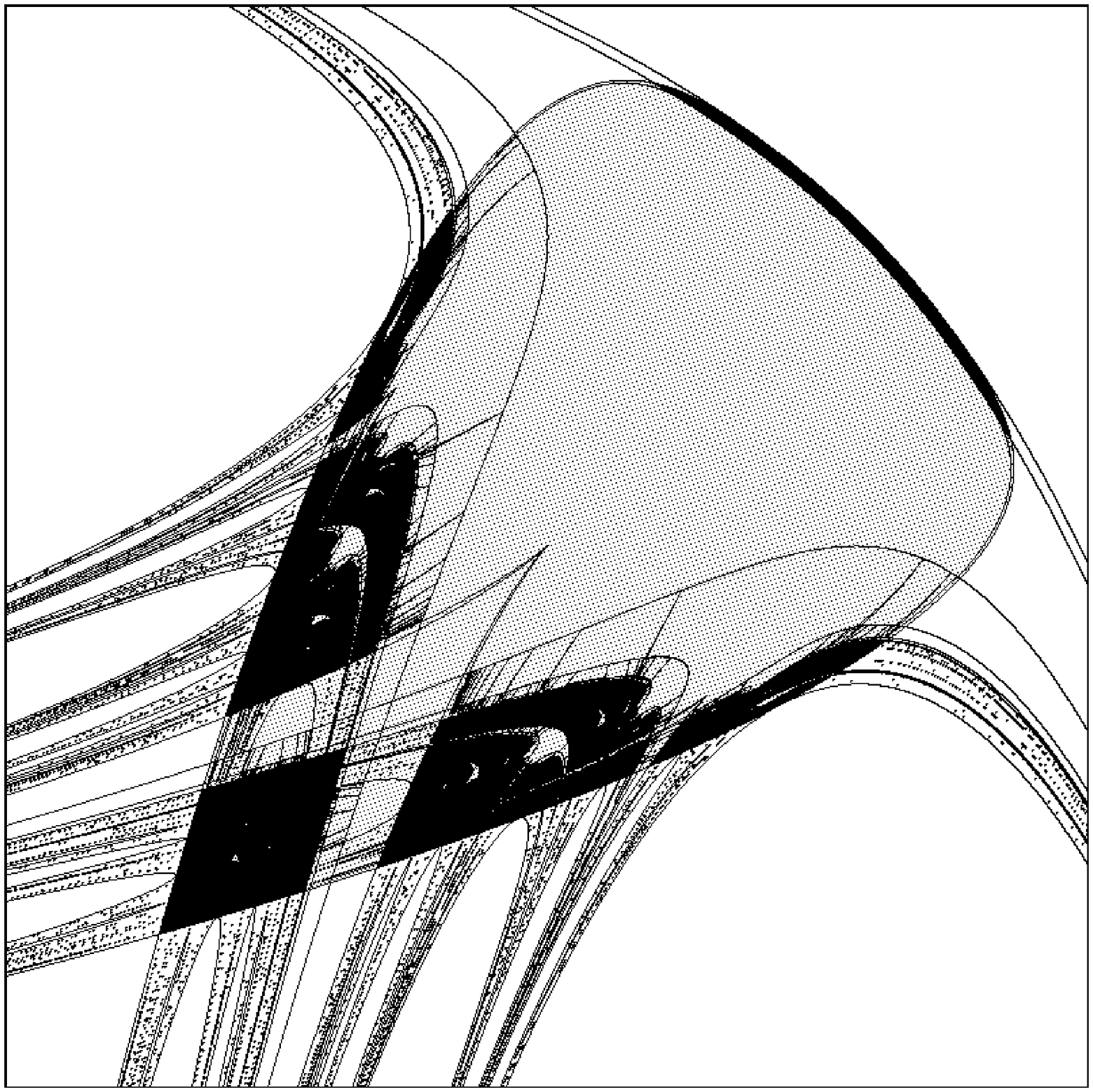} {800} {800} {225}}
\hbox{\RasterBox {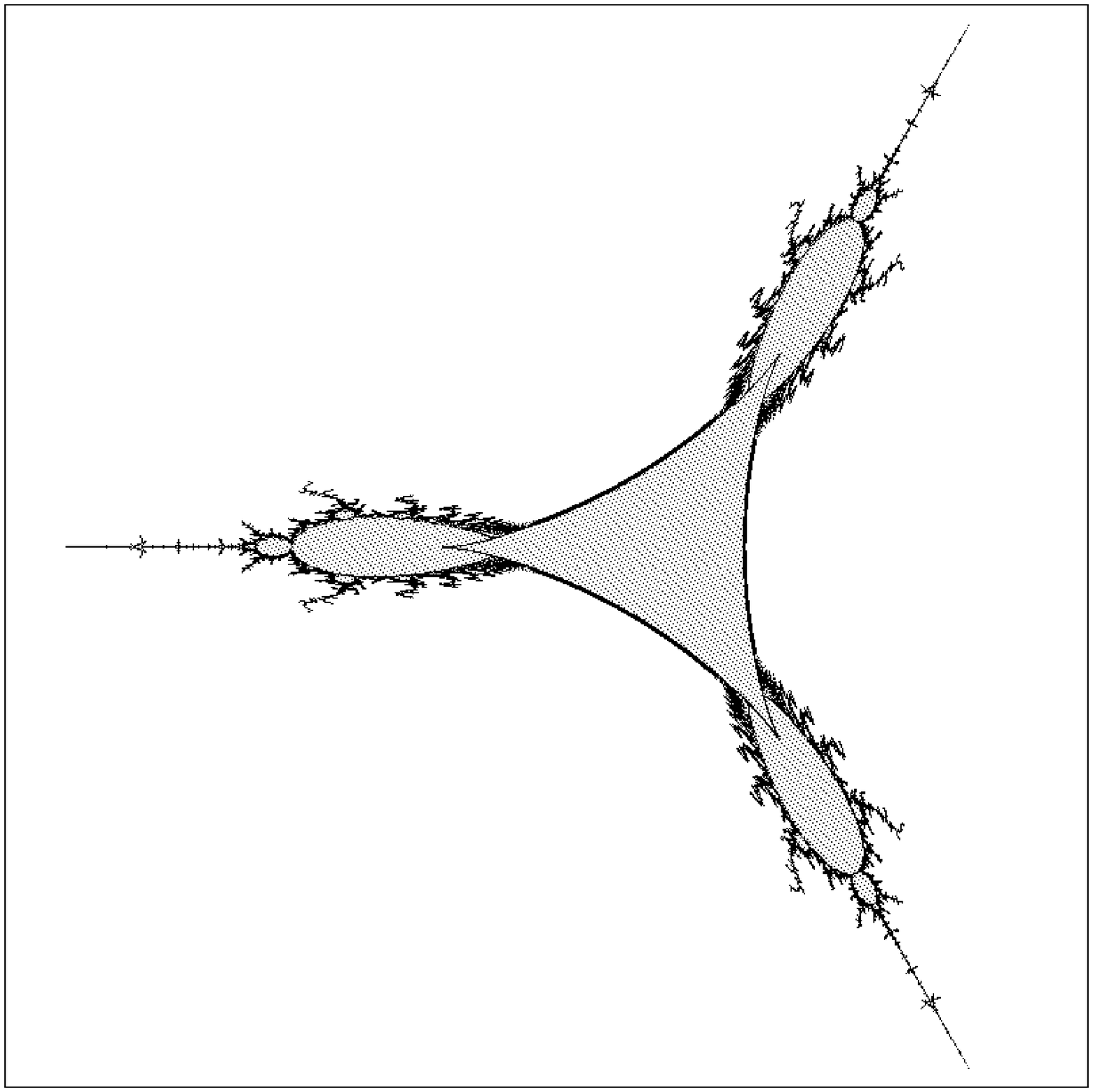} {800} {800} {225}}
\hfil}
\vskip -.1in
$$	(1,1)^+\qquad\qquad\qquad\qquad(1,1)^- $$
\hbox to \hsize{ 
\hfil  
\hbox{\RasterBox {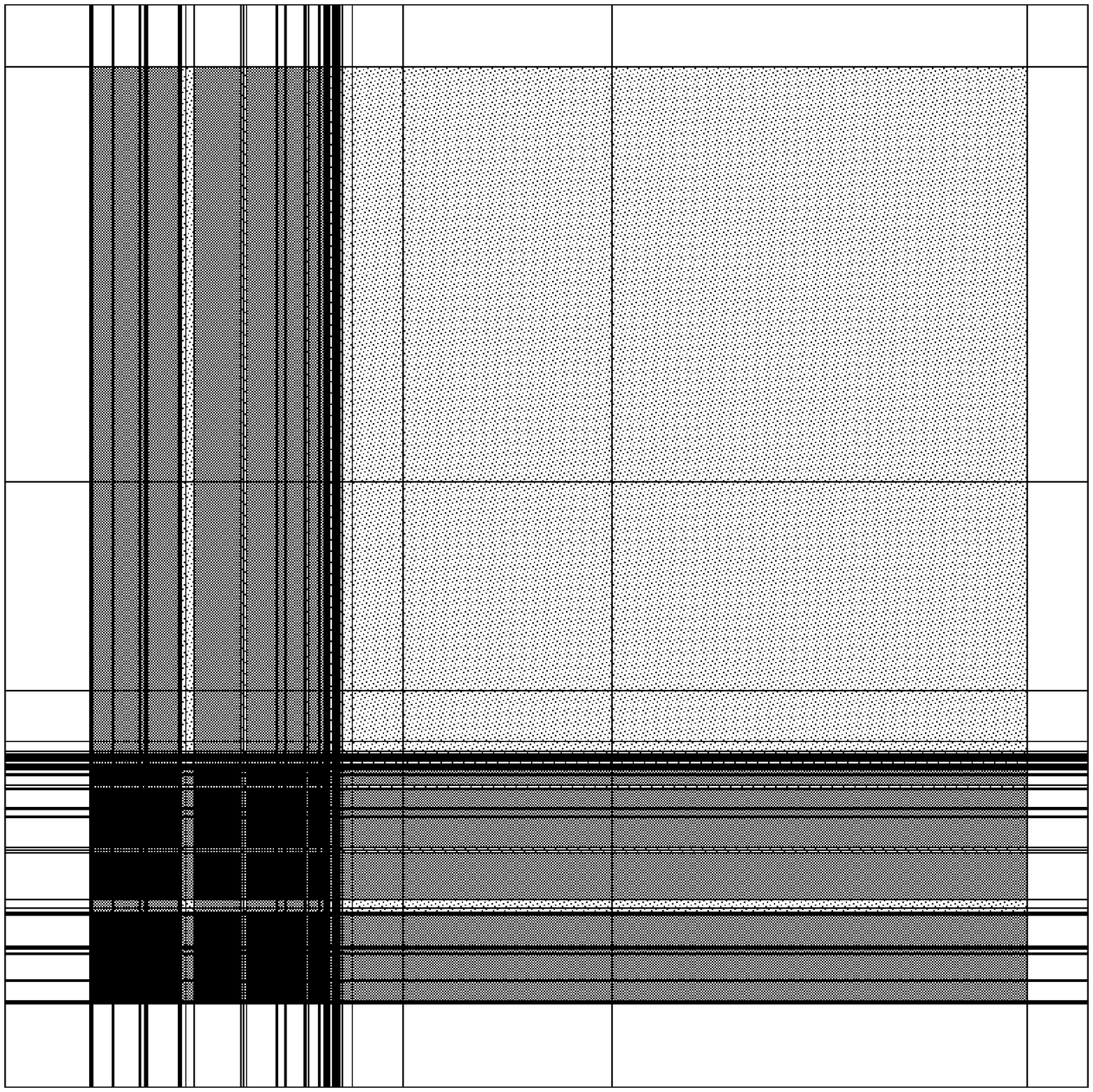} {768} {768} {225}}
\hbox{\RasterBox {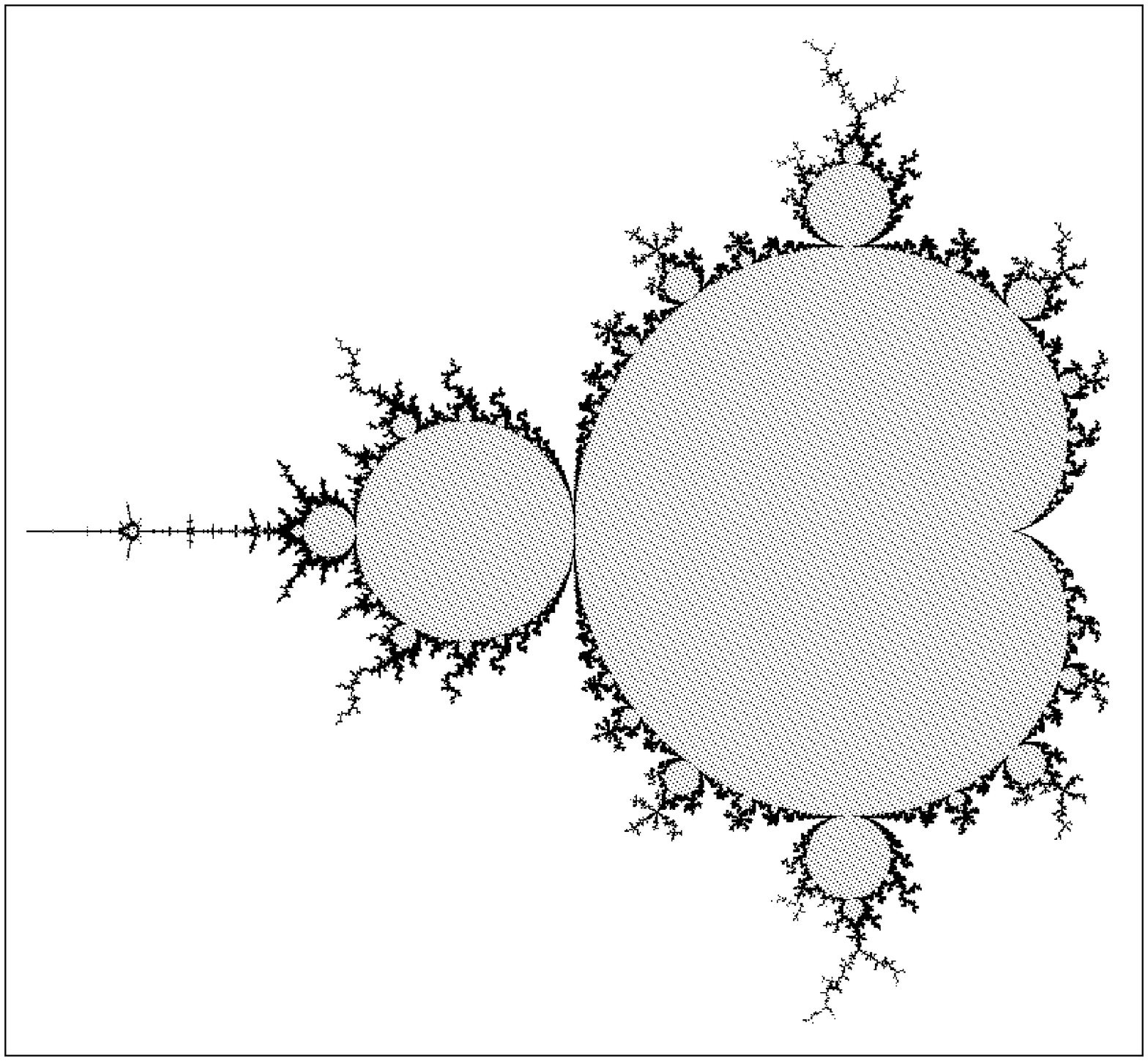} {832} {768} {225}}
\hfil}
\vskip -.1in
$$	(1)+(1)^{\;+}\qquad\qquad\qquad\qquad (1)+(1)^{\;-} $$
\centerline{\figf
Figure 4. Connectedness loci for four real forms of weight two.}
\endinsert

For each reduced mapping schema \[S\] of weight two we can make a corresponding
computation. In each case it turns out that there are exactly two real forms,
which it will be convenient to denote by \[S^+\] and \[S^-\].
(Compare Figures 3-5.) Here \[S^+\]
corresponds to the space made up out of monic real maps. For \[S^-\], we
need different descriptions in the various cases, as follows:\smallskip

For \[S^-=(1,1)^-\] and \[S^-=(1)+(1)^-\] we must consider monic centered
maps from \[|S|\times\C\] to itself which
commute with the antiholomorphic involution \[(v\,,\,z)\mapsto(3-v\,,\,\bar z)\].
Such maps have the form
$$\eqalign{(v_1\,,\,z)\mapsto(v_2\,,\,z^2+c)\;,\quad(v_2\,,\,z)\mapsto(v_1\,,\,
z^2+\bar c)\quad &{\rm for}\quad S=(1,1)\cr
(v_1\,,\,z)\mapsto(v_1\,,\,z^2+c)\;,\quad(v_2\,,\,z)
\mapsto (v_2\,,\, z^2+\bar c)\quad &{\rm for}\quad S=(1)+(1)\;.\cr}$$
The corresponding connectedness loci are shown in Figure 4. In the first case
this locus, known as the ``{\bit tricorn\/}'', has a
symmetry group \[G_\R((1,1)^-)\cong G(1,1)\] which is
non-abelian of order 6, as is visible in the Figure. (Compare [Wi].)
For \[S^-=(1)+(1)^-\], the associated connectedness locus is just a copy of the
Mandelbrot set, with the involution \[c\mapsto \bar c\] as unique real
automorphism.

\midinsert
\insertRaster 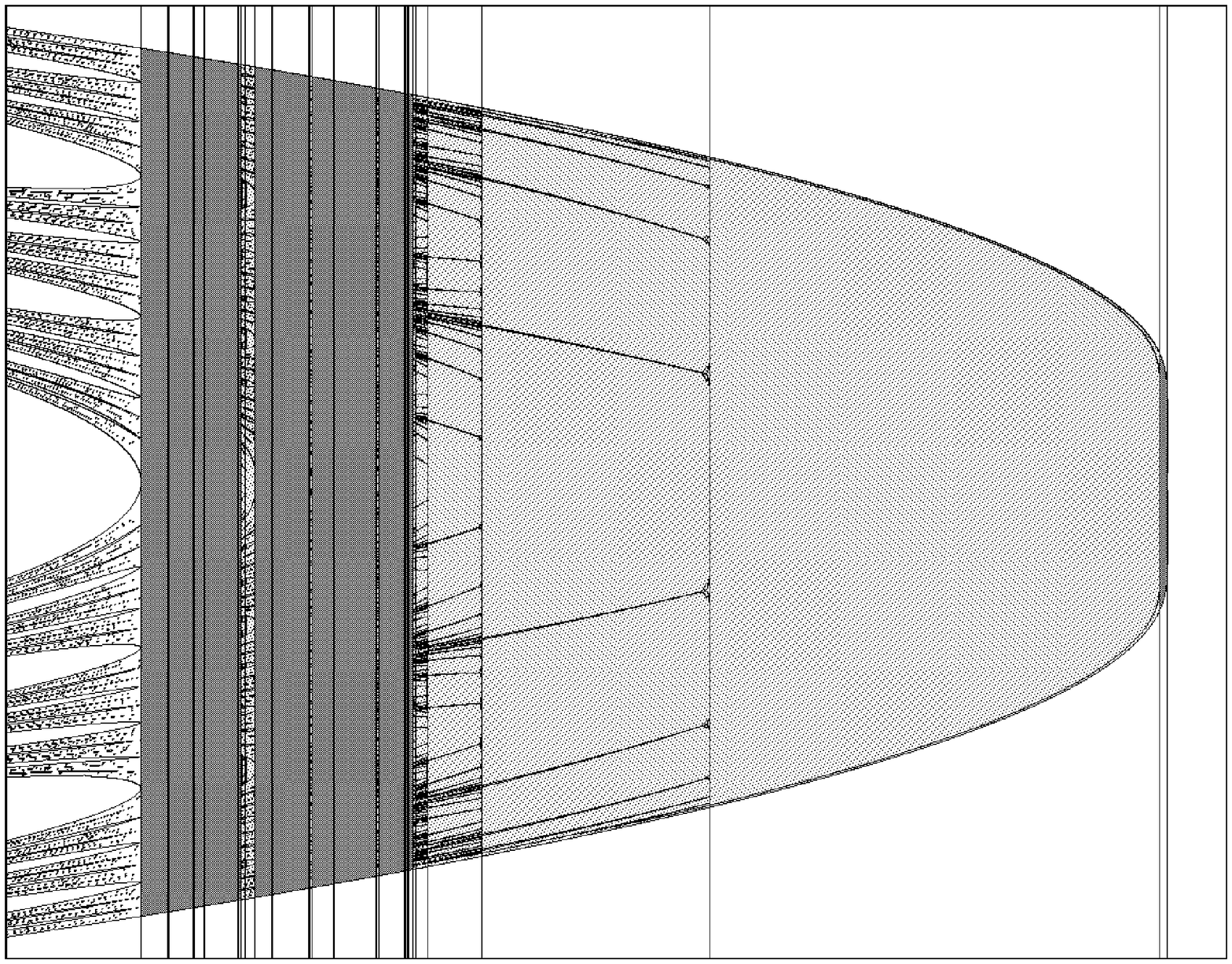 pixels 1152 by 900 scaled 225
$$	(\{1\}1)^\pm $$
\centerline{\figf
Figure 5. Connectedness loci for the remaining two real forms of weight two.}
\vfil
\endinsert

For \[S=(\{1\}1)\], the situation is more confusing, since the connectedness
loci for the two real forms \[(\{1\}1)^\pm\] are indistinguishable from
each other. We must consider real maps of the form
$$	(v_1\,,\,x)\mapsto (v_1\,,\, x^2+r_1)\;,\quad (v_2\,,\,x)\mapsto
(v_1\,,\, \pm x^2+r_2)\;, $$
where the choice of sign has very little effect on the dynamics. In this case the
group\break
\[\bar G_\R((\{1\}1)^\pm)\] of dynamic automorphisms is trivial, even though
there is an evident set theoretic
automorphism \[(r_1\,,\,r_2)\mapsto (r_1\,,\,-r_2)\] of the connectedness
locus.\smallskip

{\bf Remark 6.3.} In analogy with 2.9, there is a partial ordering
$$	(2)^\pm\quad\succ\quad\left\{\eqalign{(1,1)^- &\quad\succ\qquad
	 (1)+(1)^-\cr
	(1,1)^+&\quad\succ\quad \left\{\eqalign{&(1)+(1)^+\cr
	 &(\{1\}1)^\pm\;\;\cr}\right .\cr
	}	\right . $$
for the collection of real forms of total weight two.
\bigskip
Our main result in the real case can be stated as follows.

{\QP {\bf Theorem 6.4.} {\it Every hyperbolic component in a real
connectedness locus of weight \[w\] is a topological \[w$-cell with a unique ``center'' point,
and is real analytically homeomorphic to a uniquely defined principal
hyperbolic component \[H_{0\R}^S(
 \gamma)\], or to a suitably defined space \[B_R(S,\gamma)\] of Blaschke
products, under a homeomorphism which is uniquely determined up to the action
of the symmetry group \[\bar G_\R(S,\gamma)\].}\smallskip}

\noindent The proof involves going through the arguments in previous
sections, keeping track of the extra involution, and is not difficult.
\endproof
\vfil
\eject 

\centerline{\bf \S7. Polynomials with marked critical points.}
\smallskip

By a {\bit critically marked\/} polynomial map of degree \[w+1\] we will mean a
polynomial map \[f\] together with an ordered list \[(c_1\,,\,\ldots\,,\,
 c_{w})\] of its critical points. Even if \[f\] is a real
polynomial, this list must include all complex critical points, so that the
derivative \[f'(z)\] is a constant multiple of \[(z-c_1)\cdots(z-c_{w})\].
Similarly, we can define the concept of a ``critically
marked'' Blaschke product.

Branner and Hubbard have shown the utility of studying critically
marked polynomial
mappings. All of the principal results of the previous sections
extend to the marked case. For any mapping schema
\[S\] we can define a space \[\p^\sharp(S)\] of marked polynomial maps
and a space \[B^\sharp(S)\] made up out of marked Blaschke products.
Then any hyperbolic component of type \[S\] in any marked connectedness
locus \[C^\sharp(S')\subset \p^\sharp(S')\] is canonically homeomorphic
to \[B^\sharp(S)\]. The one step in this program which cause additional
difficulty is the analogue of Lemma 3.5, showing that \[B^\sharp(S)\] is
a topological cell. However, Douady has supplied the author
with a proof of the following key result.

{\QP {\bf Lemma 7.1.} {\it The space consisting of all critically
marked Blaschke
products of degree \[w+1\] which fix the points \[0\] and \[1\] is a
topological cell of dimension \[2w\].} \smallskip}

{\bf Proof.} We will make use of the standard diffeomorphism \[h:D\to\C\],
given by \[h(z)=w=z/\sqrt{1-|z|^2}\], with inverse \[z=w/\sqrt{1+|w|^2}\].
Let \[\beta:D\to D\] be a proper holomorphic map fixing \[0\] and \[1\],
with (not necessarily distinct) critical points \[c_1\,,\,\ldots\,,\,c_w\].
Consider the composition \[h\circ\beta:D\to\C\]. Pulling back the standard
conformal structure of \[\C\] under this composition, we obtain a new
Riemann surface
\[D_\beta\], which has underlying space \[D\] but is none-the-less conformally
isomorphic to \[\C\]. In fact we can choose the conformal isomorphism \[\eta
:\C\to D_\beta\] so that \[\eta(0)=0\] and \[\lim_{t\to+\infty}\eta(t)=1\].
The composition \[ h\circ\beta\circ\eta\] will then be a polynomial map
of the form
$$ z\;\mapsto\;h(\beta(\eta(z)))\;=\;\lambda\int_0 (w+1)\,
(z-\hat c_1)\cdots(z-\hat c_w)dz $$
with leading coefficient \[\lambda>0\], where \[\eta(\hat c_k)=c_k\]. In
fact, after replacing the function \[\eta(z)\] by \[\eta(z/
\lambda^{1/(w+1)})\], we may assume that \[\lambda=1\]. Conversely, given
such a monic polynomial with marked critical points which fixes the origin,
we can reverse this construction, thus obtaining a corresponding Blaschke
product with marked critical points. This proves that the required space
of Blaschke products is diffeomorphic to the space of \[w$-tuples
\[(\hat c_1\,,\,\ldots\,,\,c_w)\in\C^w\].\endproof

Using 3.4, it follows easily that the space of critically centered Blaschke
products
with marked critical points and with \[\beta(1)=1\] is also a topological cell.
The proof that the model space \[B^\sharp(S)\] is a topological cell is
now straightforward.\medskip

For the case of marked polynomial maps, we must expect a group of
symmetries which is much larger than \[G(S)\], and many more real forms
than in \S6. The basic definitions are as
follows. For any mapping schema \[S\], let \[\Pi(S)\] be the cartesian
product, over all vertices \[v\], of the symmetric group of permutations
of the numbers \[\{1,\,2,\,\ldots,\,w(v)\}\] which are used to index
the critical points in \[W_v\]. Then the affine
space \[\p^\sharp(S)\] of marked centered polynomial maps from \[|S|\times\C\]
to itself can be defined as a normal branched
covering space of \[\p^S\], with \[\Pi(S)\] acting as the group of deck
transformations. Similarly, the spaces \[C^\sharp(S)\supset H^\sharp_0(S)\] and \[B^\sharp(S)\]
are \[\Pi(S)$-branched coverings of \[\cl^S\supset H_0^S\]
and of \[B(S)\]. The full group \[G^\sharp(S)\] of symmetries of \[\p^\sharp(S)\] splits
as a semidirect product
$$	1 \to \Pi(S) \to G^\sharp(S) \mixarr G(S) \to 1	$$
or \[\;\; 1 \to N(S)\times \Pi(S) \to G^\sharp(S) \mixarr \Aut(S) \to 1\].
(Compare \S2.8.)\smallskip

Rather than working out the full details of this, let me simply describe the
simplest case \[S=(w)\]. If \[a_1\,,\,\ldots\,,\,a_{w}\] are arbitrary
complex numbers with barycenter\break \[a=(a_1+\cdots+a_{w})/w\], then it
may be convenient to use the notation
$$	\;z \mapsto f(z)=f(a_1\,,\,\ldots\,,\, a_{w} \,;\,z)	$$
for the unique marked centered monic map of degree \[d=w+1\] with critical
points
$$	a_1-a\,,\;\; a_2-a\,,\;\;\ldots\,,\;\;a_{w}-a	$$
and with
\[f(0)=a\]. Here the symmetric group \[\Pi(w)\] of order \[w!\] acts by
permuting the \[a_j\], and the group \[N(w)\] of \[w$-th roots of
unity \[\eta\] acts by the rule\break
\[f(z) \mapsto \eta f(z/\eta)\], or equivalently
$$	f(a_1\,,\,\ldots\,,\, a_{w} \,;\,z)\;\; \mapsto\;\; f(\eta a_1\,
 ,\,\ldots\,,\,\eta a_{w}\,;\,z)\,.	$$
The automorphism group \[\Aut(w)\] is of course trivial, so that the full
symmetry group \[G^\sharp(w)\] is just the cartesian product \[N(w)\times
 \Pi(w)\] of order \[w\,w!\].\smallskip

In order to study real forms, we extend \[G^\sharp(w)\] by the antiholomorphic
involution \[\gamma_0\], which commutes with \[\Pi(w)\],
and look for conjugacy classes of antiholomorphic
involutions in the resulting group \[\hat G^\sharp(w)\]. For \[d\] even, it turns out that there
are \[d/2\] distinct real forms, corresponding to polynomial maps with either
\[1\,,\,3\,,\,\ldots\,,\; {\rm or~}\;w\] real critical points.
For \[d\] odd, there are \[d+1\] distinct real forms. As an example, for
\[d=3\] there are two marked real forms \[\;x \mapsto \pm (x^3-3a^2x)+b\;
 \] with real critical points \[(a,-a)\] and two marked real forms
\[\;x \mapsto \pm (x^3+3a^2x)+b\;\] with imaginary critical points \[(ia,
 -ia)\]. (Compare [Branner and Hubbard].) Further
details will be left to the reader.

\end

\def\endofproof{ \smallskip \hskip 4in $\hbox{ }$ }
\input psfig
\def\QP{\narrower\medskip\noindent}

\pageno=29
\centerline{\bf Appendix}
\centerline{\bf Realizing Reduced Schemata.}
\bigskip
\centerline{by Alfredo Poirier}
\bigskip
\bigskip
We fix a reduced schema $\bar S=(|\bar S|,F,w)$. 
The main purpose of this appendix is to show that there exist a 
Postcritically Finite Polynomial $f$ of degree $d(\bar S)=w({\bar S})+1$, 
whose associated reduced schema $\bar S(f)$ realizes $\bar S$: 
\bigskip
{\bf Theorem.} 
{\it Every reduced schema $\bar S=(|\bar S|,F,w)$ can be realized.} 
\medskip
We start by constructing a non reduced schema $S$ which has $\bar S$ as associated reduced schema 
(compare Figures A.1 and A.2). 
Let $V$ be the set consisting of two copies of $|\bar S|$. 
Thus to any $v \in |\bar S|$, 
we associate $v,v' \in V$. 
We set $w_S(v)=w(v)$, 
and because we do not want to increase the degree we define $w_S(v')=0$. 
Next, 
we define the function $F_S:V \to V$ as follows. 
Let $F_S(v)=v'$ and $F_S(v')=F(v)$. 
We will realize this mapping schema $S=(|S|=V,F_S,w_S)$ 
by constructing an expanding Hubbard tree 
(for the definitions and main results in abstract Hubbard Trees we refer to [P1] or [P2]).  
\smallskip
In other words we have only inserted between every pair of (critical) vertices a third non-critical one. 
Clearly when restricted to the set $|\bar S|$, 
$F_S$ satisfies $F_S^{\circ 2}=F$. 
The next two lemmas follow directly from the above construction. 
\medskip
{\bf Lemma 1.}
{\it The mapping schema $S=(|S|,F_S,w_S)$ has reduced schema $\bar S$.} 
\endofproof
\medskip
{\bf Lemma 2.}
{\it Let $v' \in V$ be a non critical vertex. 
If $F_S(\omega)=v'$ then $v=\omega$.}
\endofproof
\smallskip
In other words every non critical vertex has exactly one preimage. 
\smallskip
\centerline{\psfig{figure=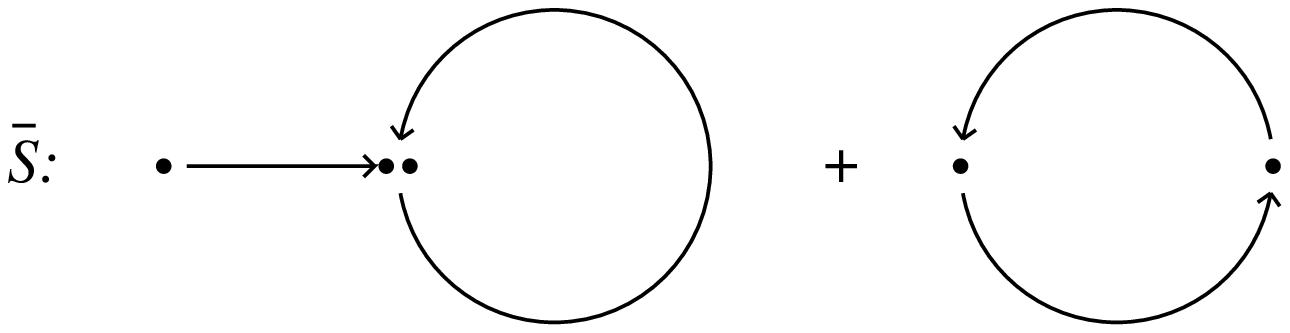,height=1.5in}}
\centerline{\figf Figure A.1: A reduced schema.}
\centerline{\psfig{figure=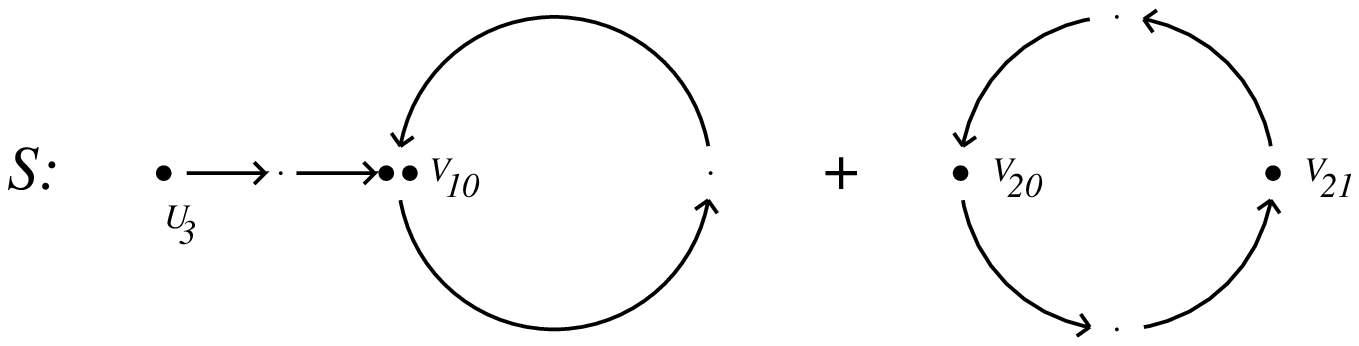,height=1.5in}}
\centerline{\figf Figure A.2: The associated non-reduced mapping schema.}
\vfill\eject
We begin the proof that the mapping schema $S$ can be realized,  
with a preliminary construction. 
Let ${\cal C}_i$ be the cycle $v_{i0} \mapsto v'_{i0} \mapsto v_{i1} \mapsto \dots v_{in}=v_{i0}$. 
For this cycle we will construct a (star-shaped) expanding Hubbard Tree $H({\cal C}_i)$ as follows. 
Join the vertices $v_{ij}$ and $v'_{ij}$ to a new vertex $p_i=p_{{\cal C}_i}$. 
We define the angle between consecutive edges at $p_i$ as $1/2n$. 
The new vertex $p_i$ will be by definition fixed of degree $d(p_i)=1$. 
At all other vertices the dynamics and degree is that induced by $S$ 
(note that by definition $d(v)=w(v)+1$). 
As there is only one periodic point which does not belong to a critical cycle, 
this is trivially an expanding abstract Hubbard tree.  
\medskip 
Next, from each critical cycle ${\cal C}_i$ 
we choose a critical vertex $v_i=v_{i0}$ and form the collection 
$\{v_1,\dots,v_m\}$; 
where $m \ge 1$ is the number of components of $S$. 
Let $u_{m+1},\dots,u_{m+r}$ be the critical vertices which do not belong to critical cycles. 
We construct an expanding abstract Hubbard tree which realizes $S$ as follows. 
Join the vertices $v_{i}$ and $u_j$ to a new vertex $q$, 
by segments $\ell_i$ and $\ell_j$. 
At $q$ these segments should form angles of non trivial multiples of $1/(r+m)$.  
Next, for every vertex $u_j$, 
we join $u_j$ and $u'_j$ by an edge which forms an angle of $1/d(u_j)$ with $\ell_j$. 
Also for every cycle ${\cal C}_i$ we paste the tree $H({\cal C}_i)$ at $v_i$ 
forming an angle of $1/d(v_i)$ with $\ell_i$. 
Because of lemma 2 this graph is a topological (angled) tree. 
If we define $q$ as a fixed vertex of degree $1$, 
then with the induced dynamics from $S$ and all $H({\cal C}_i)$ this tree is expanding 
(compare Figure A.3). 
Now, as every expanding Hubbard Tree can be realized by a Postcritically Finite polynomial (see [P1] or [P2]), 
the result follows. 
\endofproof 
\medskip
\bigskip
\centerline{\psfig{figure=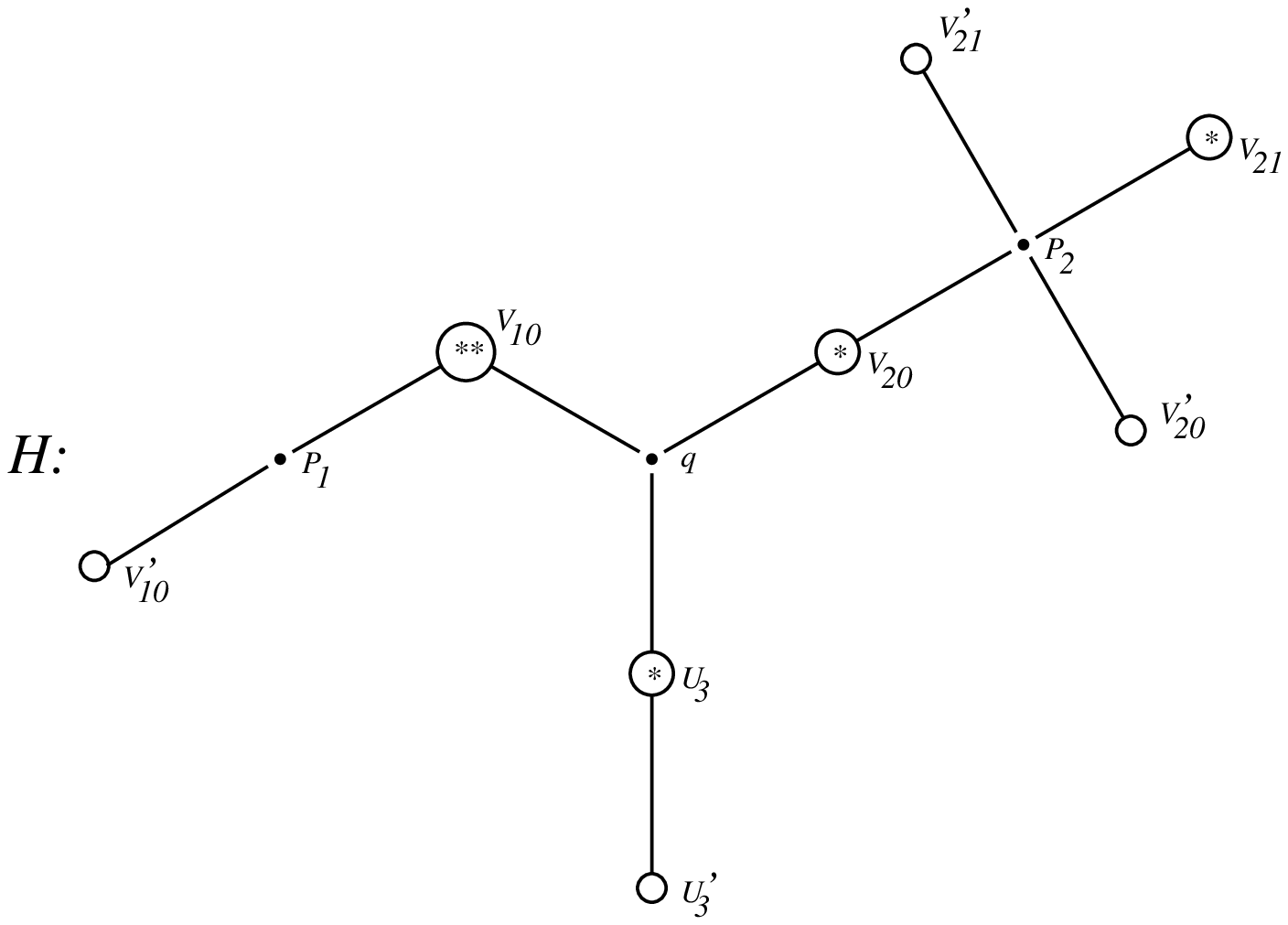,height=3in}}
\bigskip
{\QP\figf Figure A.3: An associated abstract Hubbard Tree which realizes this mapping schema.
In this figure $*$ represents a critical point 
($**$ a double critical point). 
Dots correspond to points in the Julia set, and
circles correspond to centers of Fatou components. 
As there is no edge whose endpoints correspond to Julia set vertices, this tree is 
(trivially) by definition expanding.\medskip} 
\end

\def\endproof{  \rlap{$\sqcup$}$\sqcap$ \smallskip}
\def\C{{\bf C}}
\def\QP{\narrower\medskip\noindent}
\def\[{$\,}
\def\]{\,$}
 at 12truept
\pageno=31

\def\ref{\smallskip\hangindent=1pc \hangafter=1 \noindent}
\centerline {\bf References.}   \medskip

\ref [AB] L. Ahlfors and L. Bers, The Riemann mapping theorem for variable
metrics, Annals of Math. {\bf 72} (1960), 385-404.

\ref [Bl] P. Blanchard, Complex analytic dynamics on the Riemann sphere, Bull.
Am. Math. Soc. {\bf11} (1984), 85-141.

\ref [BDK] P. Blanchard, R. Devaney and L. Keen, The dynamics of complex
polynomials and automorphisms of the shift, to appear.
 
\ref [Br] B. Branner, The parameter space for cubic polynomials, pp. 169-179
of ``Chaotic Dynamics and Fractals'', edit. Barnsley and Demko, Acad. Press
1986.
  
\ref [BD] B. Branner and A. Douady, Surgery on complex polynomials,
pp. 11-72 of ``Holomorphic Dynamics, Mexico 1986'', ed. Gomez-Mont et al.,
Lect. Notes Math. {\bf 1345}, Springer 1988.
   
\ref [BH] B. Branner and J. H. Hubbard, The iteration of cubic polynomials,
Part 1: The global topology of parameter space, Acta Math. {\bf 160} (1988)
143-206; Part 2: to appear.

\ref [BM] R. Brooks and P. Matelski, The dynamics of 2-generator subgroups
of \[PSL(2,\C)\], pp. 65-71 of ``Riemann Surfaces and Related Topics''
(Proceedings 1978 Stony Brook Conference), edit. Kra and Maskit,
Ann. Math. Stud. {\bf 97} Princeton U. Press 1981.
    
\ref [D1] A. Douady, Syst\'emes dynamiques holomorphes, S\'em. Bourbaki 599,
Ast\'erisque {\bf 105-106} (1983), 39-63.

\ref [D2] A. Douady, L'\'etude dynamique des polyn\^omes quadratiques complexes
et ses reinvestissements, Soc. Math. France, Assembl. Gen. Jan. 1985, 21-42.

\ref [DE] A. Douady and C. Earle, Conformally natural extension of homeomorphisms
of the circle, Acta math. 157 (1986), 23-48.

\ref [DH1] A. Douady and J. H. Hubbard, \'Etude dynamique des polyn\^omes
complexes I \& II, Publ. Math. Orsay (1984-85).
     
\ref [DH2] A. Douady and J. H. Hubbard, On the dynamics of polynomial-like mappings,
Ann. Sci. Ec. Norm. Sup. (Paris) {\bf 18} (1985), 287-343.

\ref [La] P. Lavaurs, Syst\`emes dynamiques holomorphes: explosion de points
p\'eriodiques para\-boliques, Thesis, Univ. Paris-Sud Orsay 1989.

\ref [LV]O. Lehto and K. I. Virtanen, Quasiconformal Mappings in the Plane,
Springer 1973.

\ref [Ly1] M. Lyubich, The dynamics of rational transforms: the topological picture,
Russian Math. Surveys {\bf 41:4} (1986), 43-117.

\ref [Ly2] M. Lyubich, Some typical properties of the dynamics of rational
maps, Russian Math. Surveys {\bf 38} (1983) 154-155.

\ref [Ly3] M. Lyubich, An analysis of the stability of the dynamics of rational
functions, Funk. Anal. i. Pril. {\bf 42} 1984), 72-91; Selecta Math. Sovietica
{\bf 9} (1990) 69-90.

\ref [MSS]
R. Ma\~n\'e, P. Sad and D. Sullivan, On the dynamics of rational maps,
Ann. Sci. \'Ec. Norm. Sup. Paris, {\bf 16} (1983), 193-217.

\ref [McM] C. McMullen, Automorphisms of rational maps, pp. 31-60 of
``Holomorphic Functions and Moduli I'', ed. Drasin, Earle, Gehring, Kra \&
Marden; MSRI Publ.$\,10$, Springer 1988.

\ref [M1] J. Milnor, Self-similarity and hairiness in the Mandelbrot set,
pp. 211-257 of ``Computers in Geometry and Topology'', edit. Tangora,
Lect. Notes Pure Appl. Math. {\bf 114}, Dekker 1989.

\ref [M2] J. Milnor, Dynamics in one complex variable: Introductory lectures,
Stony Brook IMS Preprint 1990/5.

\ref [M3] J. Milnor, Remarks on iterated cubic maps, Stony Brook IMS
Preprint 1990/6 (to appear in Experimental Mathematics).

\ref [M4] J. Milnor, Remarks on quadratic rational maps,
in preparation.

\ref [P1] A. Poirier, On the realization of fixed point portraits,
Stony Brook IMS Preprint 1991/20.

\ref [P2] A. Poirier, On postcritically finite polynomials,
Thesis, Stony Brook, in preparation.

\ref [P3] A. Poirier, Hubbard forests, in preparation.

\ref [R1] M. Rees, Components of degree two hyperbolic rational maps,
Invent. Math. {\bf 100} 1990, 357-382.

\ref [R2] M. Rees, { A partial description of parameter space of rational
maps of degree two}, Part I: preprint, Univ. Liverpool 1990;
and Part II: Stony Brook IMS preprints 1991/4.

\ref [Sh1] M. Shishikura, Surgery of complex analytic dynamical systems,
pp. 93-105 of ``Dynamical Systems and Nonlinear Oscillations'', ed.
G. Ikegami, World Scient. Series in Dyn. Sys. {\bf 1}, Singapore 1986.

\ref [Sh2] M. Shishikura, On the quasiconformal surgery of rational functions,
Ann. Sci. \'Ec. Norm. Sup. Paris, {\bf 20} (1987) 1-29.

\ref [Su1] D. Sullivan, Conformal dynamical systems, pp. 725-752 of ``Geometric
Dynamics'', ed. J. Palis, Lect. Notes Math. {\bf 1007}, Springer 1983.

\ref [Su2] D. Sullivan, Quasiconformal homeomorphisms in dynamics, topology,
and geometry,\break preprint 1986.

\ref [Su3] D. Sullivan, Quasiconformal homeomorphisms and dynamics, I,
Annals of Math. {\bf 122} (1985) 401-418; II, Acta Math. {\bf 155} (1985)
243-260; III, to appear.

\ref [Su4]
D. Sullivan, Expanding endomorphisms in the complex plane, to appear.

\ref [Wi]
R. Winters, Bifurcations in families of antiholomorphic and biquadratic
maps, Thesis, Boston University 1990.


\end